\documentclass[journal]{IEEEtran}

\usepackage{ifpdf}
\usepackage{cite}
\usepackage{amsmath, nccmath}
\usepackage{url}
\usepackage{graphicx}
\usepackage{mathtools}

\renewcommand*{\arraystretch}{1.2}
\graphicspath{{graphics/}}

\begin{document}

\title{Ready-to-Use Unbiased Estimators for Multivariate Cumulants Including One That Outperforms $\overline{x^3}$}
\author{Fabian~Schefczik
        and~Daniel~H\"agele
\thanks{F. Schefczik and D. H\"agele are with the Institute for Experimental Physics VI, Ruhr-University Bochum,
Germany, e-mail: (see http://www.optics.rub.de).}
\thanks{submitted April 2019.}}
\markboth{IEEE TRANSACTIONS ON SIGNAL PROCESSING, submitted April 2019}%
{Shell \MakeLowercase{\textit{et al.}}: Estimators for Higher-Order Multivariate Cumulants}
%


\maketitle

\begin{abstract}
We present multivariate unbiased estimators for second, third, and fourth order cumulants $C_2(x,y)$, $C_3(x,y,z)$, and $C_4(x,y,z,w)$.  Many relevant new estimators are derived for cases where some variables are average-free or pairs of variables have a vanishing second order cumulant. 
The well-know Fisher k-statistics is recovered for the single variable case.
The variances of several estimators are explicitly given in terms of higher order cumulants and discussed with respect to random processes that are predominately Gaussian. We surprisingly find that the frequently used third order estimator $\overline{x^3}$ for $C_3(x,x,x)$ of a process $x$ with zero average is outperformed by alternative estimators. The new (Gauss-optimal) estimator $\overline{x^3} - 3 \overline{x^2}\overline{x}(m-1)/(m+1)$ improves the variance by a factor of up to $5/2$. Similarly, the estimator 
$\overline{x^2 z}$ for $C_3(x,x,z)$ can be replaced by another Gauss-optimal estimator. The known estimator $\overline{xyz}$ for $C_3(x,y,z)$ as well as
previously known estimators for $C_2$ and $C_4$ of one average-free variable are shown to be Gauss-optimal.
As a side result of our work we present two simple recursive formulas for finding multivariate cumulants from moments and vice versa.
\end{abstract}

\begin{IEEEkeywords}
bias, consistency, cumulant, estimation,
estimator, higher moments, higher order statistics.
\end{IEEEkeywords}
%
\IEEEpeerreviewmaketitle

\section{Introduction}
\IEEEPARstart{H}{igher order cumulants} find various applications in signal processing for the investigation of higher order statistics \cite{comonSP1994}. Methods for blind source separation heavily rely on higher order cumulants of stochastic vectors, i.e. on cumulants of more than one variable. Polyspectra for higher order harmonic analysis are another example of a concept that is based on higher order multivariate cumulants \cite{brillingerAMS1965}. Third and fourth order polyspectra recently found application in physics for the investigation of continuous quantum measurements \cite{hagelePRB2018}. 

The calculation of a cumulant like, e.g., the variance $C_2(x,x) = \langle x^2 \rangle - \langle x \rangle^2$ from a finite sample of $m$ data $x_j$ requires the choice of a suitable estimator. The formula 
\begin{eqnarray}
   k_2 & = & \frac{1}{m-1} \sum_{j = 1}^m \left(x_j - \left(\frac{1}{m}\sum_{j'=1}^m x_{j'}\right)\right)^2 \nonumber \\
    & = & \frac{m}{m-1}\overline{(x-\overline{x})^2} \label{bessel} 
\end{eqnarray}
with the famous Bessel correction $m-1$ in the denominator is known to be an unbiased estimator for a non-average free independent and identically distributed (i.i.d.) process. 
In the second line, we denoted the average of $m$-samples by $\overline{\cdots}$. This notation allows for more compact equations and will be used throughout the article.
It can be shown that $C_2(x,x) = \langle k_2 \rangle$ while the statistical expectation of the so-called classical or natural estimator $\overline{(x-\overline{x})^2}$ deviates from $C_2(x,x)$ by an error that is on the order of $1/m$.
In 1928 R. A. Fisher gave  corresponding explicit formulas for $k_1$ to $k_6$ and a general recipe for obtaining even higher order unbiased estimators for cumulants of one or more variables \cite{FisherPLMS1928}. His formulas $k_j$ are today known as k-statistics 
and find frequent application in signal processing. Nardo {\it et al.} gave in 2009 a general framework for deriving multivariate
 k-statistics based on the mathematical concept of umbra calculus \cite{DiNardoBernoulli2008,rotaSIAM1994}. 
 They, however, give only one ready-to-use multivariate estimator for the cumulant $C_3(x,x,z)$ with $x$ and $z$ being not average-free. 
Mansour as well as Blagouchine and Moreau recently derived and discussed an unbiased estimator for the fourth order cumulant of an average-free variable where $\langle x \rangle = 0$ was exploited to yield a simpler estimator as compared to the general $k_4$ \cite{mansourCONF1998,BlagouchineIEEE2009}.
Despite the obvious demand of a generalized version of k-statistics for the multivariate case or for cases with additional knowledge on the random variables, we are aware of only a small number of  special cases in the literature (compare text below Table \ref{OverviewEstimator}).
Here we present a collection of unbiased estimators for multivariate cumulants up to fourth order 
 including cases where one or more variables are average free or the covariance $C_2$ of a pair of different variables is known to vanish.
These estimators contain less terms than the general (or full) estimators and will therefore be referred to as {\it reduced estimators}. We recover the known result that 
the reduced estimator of $C_4(x,x,x,x)$ exhibits an improved variance for a Gauss distributed variable compared to $k_4$, i.e. the signal to noise improves upon using the reduced estimator  \cite{BlagouchineIEEE2009}. Much to our surprise, we find that the corresponding reduced estimator $\overline{x^3}$ for $C_3(x,x,x)$ is in fact worse than $k_3 = m^2\overline{(x - \overline{x})^3}/((m-1)(m-2))$. Below, a new (Gauss-) optimal estimator is derived that even surpasses the performance of $k_3$. We expect that the new estimators will soon replace biased or non-Gauss-optimal estimators that have been used in literature for the  lack of alternatives. 

The paper is organized as follows. We first give a short review on multivariate higher order cumulants and some of their properties in Section \ref{sec:cum}. In Sections \ref{sec:c2}, \ref{sec:c3}, \ref{sec:c4}, and \ref{sec:c4c} we derive unbiased estimators for second to fourth order cumulants for random variables with various known properties like zero average. In Section \ref{sec:var} we discuss the variance of several estimators and are led to the question of estimators that are optimal under the condition of random variables that are predominantly Gaussian. In Section \ref{sec:go} we derive several Gauss-optimal estimators. Their usefulness is illustrated by
a numerical example in Section \ref{sec:num} for an average free random variable with slightly asymmetric distribution. 
Owing to the considerable number of different estimators, we introduce here a new nomenclature for them based on 
conditions for their application. Table \ref{OverviewEstimator} gives an overview of all estimators including references for those
 previously known. 

\section{Cumulants}
\label{sec:cum}
The $n$th order cumulants of a stochastic vector $\vec{x} = (x_1, x_2, \cdots )^{\rm T}$ can be defined by a generating function
\begin{equation}
	K_{\vec{x}} (\vec{k}) = \ln\langle \exp \left(\vec{k} \cdot \vec{x}\right) \rangle
\end{equation}
and its derivatives at $\vec{k} = 0$
\begin{equation}
	C_n\left(x_1,...,x_n\right) = \frac{\partial^n}{\partial k_1...\partial k_n} 	K_{\vec{x}} (\vec{k}) \Big|_{\vec{k} = 0}.
\end{equation}
Here we used angle brackets $\langle \cdots \rangle$ to denote the expected statistical value, instead of $E[\cdots]$ 
to obtain slightly more compact expressions. 
The definitions $x_1=x$, $x_2=y$, $x_3=z$, and $x_4=w$ are used throughout the paper for the same reason. The four lowest-order multivariate cumulants are then explicitly given in terms of products of higher order moments \cite{gardinerBOOK2009} 
\begin{eqnarray}
	C_1(x) &=& \langle x \rangle \label{C_1} \\
	C_2(x,y) &=& \langle xy \rangle - \langle x\rangle \langle y \rangle \label{C_2}\\
\nonumber	C_3(x,y,z) &=& \langle xyz \rangle - \langle xy \rangle \langle z \rangle - \langle xz \rangle \langle y \rangle - \langle yz \rangle \langle x \rangle \\
 && + 2 \langle x \rangle \langle y \rangle \langle z \rangle \label{C_3}\\
\nonumber C_4(x,y,z,w) &=&  \langle xyzw \rangle  - \langle xyz \rangle \langle w \rangle - \langle xyw \rangle \langle z \rangle\\
\nonumber && - \langle xzw \rangle \langle y \rangle - \langle yzw \rangle \langle x \rangle - \langle xy \rangle \langle zw \rangle \\
\nonumber && - \langle xz \rangle \langle yw \rangle - \langle xw \rangle \langle yz \rangle \\
\nonumber && +2 \langle xy \rangle \langle z \rangle \langle w \rangle +2 \langle xz \rangle \langle y \rangle \langle w \rangle  \\
\nonumber && +2 \langle xw \rangle \langle y \rangle \langle z\rangle +2 \langle yz \rangle \langle x \rangle \langle w \rangle \\
\nonumber && +2 \langle yw \rangle \langle x \rangle \langle z \rangle +2 \langle zw \rangle \langle x \rangle \langle y \rangle \\
&& - 6 \langle x \rangle \langle y \rangle \langle z \rangle \langle w \rangle. \label{C_4}
\end{eqnarray}
We give in Appendix \ref{app:var} a recursive relation between cumulants and moments that can be used to obtain multivariate cumulants of any order in terms of moments and vice versa. 
Any cumulant $C_n$ of the sum of two independent stochastic vectors $\vec{x}$ and $\vec{y}$ shows the important property
\begin{equation}
	C_n(\vec{x}+\vec{y}) = C_n(\vec{x}) + C_n(\vec{y}).
\end{equation}
If $\vec{x}$ is a desired signal and $\vec{y}$ an undesired background noise (e.g. electronic noise of an amplifier), then
$C(\vec{x})$ can be determined via that above relation from $C(\vec{x}+\vec{y})$ and a separately measured background cumulant 
 $C_n(\vec{y})$. 
In contrast to cumulants, such a procedure is not possible using higher order moments as for $M(\vec{x}) = \langle x_1 x_2  ... \rangle$ 
\begin{equation}
	M_n(\vec{x}+\vec{y}) \neq M_n(\vec{x}) + M_n(\vec{y})
\end{equation}
for $n \ge 2$. This is the main reason why cumulants are so important for the evaluation of actual experiments.
 In the following we will derive estimators $c_n(\vec{x})$ of the cumulants $C_n(\vec{x})$. The estimators are functions of
  $m$ samples of $\vec{x}$ and will be constructed in a way that they fulfill  $C_n(\vec{x}) = \langle c_n(\vec{x}) \rangle$. Any such 
 an estimator is called unbiased, in contrast to a biased estimator $\tilde{c}_n(\vec{x})$ with an error $\varepsilon(m)$ where  $C_n(\vec{x}) = \langle \tilde{c}_n(\vec{x}) \rangle + \varepsilon(m)$. 
\begin{table*}
\caption{\label{OverviewEstimator}  unbiased reduced and Gauss-optimal estimators.}
\begin{tabular}{l l l l}
Conditions & $C_2(x,y)$  & $C_3(x,y,z)$ & $C_4(x,y,z,w)$ \\
\hline
& & & \\
(a) no conditions & $c_2^{\rm (a)}(x,y)$ GV & $c_3^{\rm (a)}(x,y,z)$ GV & $c_4^{\rm (a)}(x,y,z,w)$ G\\
(b) $x = y = z = w$ & $c_2^{\rm (b)}(x) = k_2$ GV & $c_3^{\rm (b)}(x) = k_3$ GV & $c_4^{\rm (b)}(x) = k_4$ GV \\
(c) $x = y = z = w$  and $\langle x \rangle = 0$ &
                  $c_2^{\rm (c)}(x) = c_2^{\rm (e)}(x,x)$ GV & $c_3^{(c)}(x)$ nG V, $c_3^{\rm (c,Go)}(x)$ GV$^*$ & $c_4^{\rm (c)}(x)$ GV \\
(d) $\langle x \rangle = \langle y \rangle = \langle z \rangle = \langle w \rangle = 0$ & $c_2^{\rm (d)}(x,y) = c_2^{\rm (e)}(x,y)$ GV & $c_3^{\rm (d)}(x,y,z)$ GV
 & $c_4^{\rm (d)}(x,y,z,w)$  \\
(e) only $\langle x \rangle = 0$ & $c_2^{\rm (e)}(x,y)$ GV & $c_3^{\rm (e)}(x,y,z)$ GV$^*$ &            $c_4^{\rm (e)}(x,y,z,w)$ \\
(f) only $\langle x \rangle = \langle y \rangle = 0$  & - & $c_3^{\rm (f)}(x,y,z)$ GV$^*$ & $c_4^{\rm (f)}(x,y,z,w)$      \\
(g) only $\langle x \rangle = \langle y \rangle = \langle z \rangle = 0$ & - & - &   $c_4^{\rm (g)}(x,y,z,w)$ \\
(h) $\langle x \rangle = \langle y \rangle = \langle z \rangle = \langle w \rangle = 0$ and $x=y$ & - & $c_3^{\rm (h)}(x,z) = c_3^{\rm (d)}(x,x,z)$ nG V, 
 & not treated \\
 & & $c_3^{\rm (h,Go)}(x,z)$ GV$^*$ & \\
(i) $C_2(x,y) = 0$ & - & $c_3^{\rm (i)}(x,y,z)$ & not treated \\
(j) $\langle x \rangle = \langle y \rangle = \langle z \rangle = \langle w \rangle = 0$  & 
- & $c_3^{\rm (j)}(x,y,z)$ & $c_4^{\rm (j)}(x,y,z,w)$ \\
and $C_2(x,y)=0$, $C_2(x,z)=0$, etc. & & & \\
(ca) $x=a$, $y = a^*$, $z=b$, $w = b^*$& - & - & $c_4^{\rm (ca)}(a,b)$ V \\
and $\langle a \rangle = \langle b \rangle = 0$, $C_2(a , b) = 0$, $C_2(a,b^*) = 0$ & & & \\
(cb) $x=a$, $y = a^*$, $z=b$, $w = b^*$& - & - & $c_4^{\rm (cb)}(a,b)$ V$^*$ \\
and $\langle a \rangle = \langle b \rangle = 0$, $C_2(a , b) = 0$, $C_2(a,b^*) \neq 0$ & & & \\
(cc) $x=a$, $y = a^*$, $z=b$, $w = b^*$& - & - & $c_4^{\rm (cc)}(a,b)$  \\
and $\langle a \rangle \neq 0$, $\langle b \rangle = 0$, $C_2(a , b) = 0$, $C_2(a,b^*) = 0$ & & & \\
(cd) $x=a$, $y = a^*$, $z=b$, $w = b^*$ & - & $c_3^{\rm (cd)}(a, b)$ GV$^*$ & -  \\
and $\langle a \rangle = 0$, $\langle b \rangle \neq 0$, $C_2(a , b) = 0$, $C_2(a,b^*) = 0$ & & & \\
& & & 
\end{tabular}
 \\
The table gives an overview on which estimator can be used to estimate a cumulant depending on conditions.
The univariate estimators $c_j^{\rm (b)}(x)$ are identical with Fisher's k-statistics $k_j$ \cite{FisherPLMS1928}. 
The estimator $c_4^{\rm (c)}$ has been introduced and treated in \cite{mansourCONF1998,BlagouchineIEEE2009} and $c_4^{\rm (ca)}(a,b)$ was introduced in
\cite{starosielecBOOK2012} and applied in \cite{starosielecRSI2010}.
 Estimators labeled with G are proved to be Gauss-optimal in the text. Estimator labeled nG are proved to be not Gauss-optimal. The general variance of the estimator is given in the text if it is labeled with V. The variance of the estimator for Gaussian variables $x$, $y$, is given if it is labeled V$^*$. All estimators are consistent (Section \ref{sec:consistency}). 
\end{table*}   
\section{Estimators for $C_2$}
\label{sec:c2}
While the unbiased estimators for $C_2$ are well known, they are derived here for didactical reasons. The same method used here  will be applied to derive estimators of $C_3$ and $C_4$ in the following sections.
We seek to write $C_2(x,y)$ as the expected statistical value of first and second order means of $m$ samples of $\vec{x}$.
The expected statistical value of the mean  of $m$ samples $x_i y_i$
exhibits a simple relation with the second order moment $\langle xy \rangle$ 
\begin{equation}
	\langle \overline{xy} \rangle = \frac{1}{m} \sum_{i}^{m} \langle x_iy_i \rangle = \langle xy \rangle,
\end{equation}
where we used an overline $\overline{\cdots}$ to denote the mean of $m$ samples. 
Consider now
\begin{equation}
	\langle \overline{x} \ \overline{y} \rangle = \frac{1}{m^2}  \sum_{i,j}^{m}\langle x_iy_j \rangle,
\end{equation}
which cannot be reduced to a single expected value. There are $m$ terms under the sum with $i = j$ and $m(m-1)$ terms with
 $i \neq j$. The first case gives rise to contributions $\langle x y \rangle$ and the second to contributions $\langle x \rangle \langle y \rangle$ (the second case requires $\vec{x}$ to be i.i.d. !) which yields
\begin{equation}
	\langle \overline{x} \ \overline{y} \rangle = \frac{1}{m^2} \left( m\langle xy \rangle + m(m-1)\langle x \rangle \langle y \rangle \right).
\end{equation}
The above relations between the sample means and the expected values can be combined into a single matrix equation 
\begin{equation}
\begin{pmatrix}
\langle \overline{xy}\rangle\\
\langle \overline{x} \ \overline{y} \rangle\\
\end{pmatrix}
=  \begin{pmatrix}
1 & 0  \\
\frac{m}{m^2} & \frac{m(m-1)}{m^2}  \\
\end{pmatrix} 
\begin{pmatrix}
\langle xy \rangle \\
\langle x \rangle \langle y \rangle\\
\end{pmatrix}. \label{2matrix}
\end{equation}
The second-order cumulant (\ref{C_2}) can be written as
\begin{equation}
C_2(x,y)
 = \begin{pmatrix}
1 & -1 
\end{pmatrix} 
\begin{pmatrix}
\langle xy \rangle \\
\langle x \rangle \langle y \rangle\\
\end{pmatrix},  \label{C2def}
\end{equation}
and expressed in terms of sample means denoting the above $2\times2$-matrix by $A_2$ as
  \begin{equation}
C_2(x,y) = \left\langle \begin{pmatrix}
1 &-1
\end{pmatrix} A_2^{-1}
\begin{pmatrix}
\overline{xy} \\
\overline{x} \ \overline{y}\\
\end{pmatrix} \right\rangle.  \label{C2def2}
\end{equation}
We can find an unbiased estimator for $C_2(x,y)$ in the angle brackets of the above equation 
\begin{equation}
	c_2^{\rm (a)}(x,y) = \frac{m}{m-1} \left(\overline{xy} - \overline{x} \ \overline{y}\right)
	\label{Variance}
\end{equation}
and for the case of just one variable
\begin{equation}
c_2^{\rm (b)}(x) = \frac{m}{m-1} \left(\overline{x^2} - \overline{x}^2\right)
\label{Variance_x}
\end{equation}
which is identical to (\ref{bessel}) and correctly exhibits the Bessel correction.  
If the processes are known to have zero means $\langle x \rangle =  0$ and/or $\langle y \rangle= 0$,
 (\ref{C2def}) simplifies to 
\begin{equation}
C_2(x,y) = \begin{pmatrix}
1 & 0
\end{pmatrix} \cdot
\begin{pmatrix}
\langle xy \rangle \\
\langle x \rangle \langle y \rangle\\
\end{pmatrix} \label{eq:parameterZero}
\end{equation}
 and a reduced estimator
\begin{equation}
c_2^{\rm (d/e)}(x,y) =  \overline{xy} 
	\label{ReducedVariance}
\end{equation}
follows without the Bessel correction.
For $x=y$ we have
\begin{equation}
c_2^{\rm (c)}(x) =  \overline{x^2}.
\end{equation}
Please note that the zero in (\ref{eq:parameterZero}) can be replaced by any parameter $\alpha$ since $\langle x \rangle = 0$. This would 
yield correct unbiased estimators of $C_2$ for any $\alpha$. We will see below that sometimes an $\alpha \neq 0$ can be found that 
yields an improved variance compared to the estimator with $\alpha = 0$.    

\section{Estimators for $C_3$}
\label{sec:c3}
Next, we derive unbiased estimators for the multivariate cumulant $C_3$ [see Eq. (\ref{C_3})]. The derivation follows the same scheme as above. 
 Considering
\begin{equation}
\langle \overline{xyz} \rangle = \frac{1}{m} \sum_{i}^{m} \langle x_i y_i z_i \rangle = \langle xyz \rangle,
\end{equation}
we establish $\overline{xyz}$ as an unbiased estimator of $\langle xyz \rangle$.
Considering
\begin{equation}
	\langle \overline{xy} \ \overline{z} \rangle  = \frac{1}{m^2} \sum_{i,j}^{m} \langle x_i y_i z_j \rangle,
\end{equation}
we find that the terms under the sum appear with the multiplicities
\begin{eqnarray}
\nonumber i=j=k \ &\textrm{$m$-times} \\
\nonumber i=j\neq k  \ &\textrm{$m(m-1)$-times} \\
\end{eqnarray}
which leads us to
\begin{equation}
m^2\langle \overline{xy} \ \overline{z} \rangle =  m(m-1) \langle xy \rangle \langle z \rangle+ m \langle xyz \rangle.
\end{equation}
Corresponding expressions hold for  $\langle \overline{xz} \ \overline{y} \rangle $ and $\langle \overline{yz} \ \overline{x} \rangle$.

Last, we treat
\begin{equation}
\langle \overline{x} \ \overline{y} \ \overline{z} \rangle  = \frac{1}{m^3} \sum_{i,j,k}^{m} \langle x_i y_j z_k \rangle
\end{equation}
where we find the multiplicities
\begin{eqnarray}
\nonumber i=j=k  \ &\textrm{$m$-times} \\
\nonumber i=j\neq k  \ &\textrm{$m(m-1)$-times}  
\nonumber i\neq j = k  \ &\textrm{$m(m-1)$-times} \\
\nonumber i=k \neq j  \ &\textrm{ $m(m-1)$-times} \\
\nonumber i\neq j\neq k \ &\textrm{$m(m-1)(m-2)$-times,} 
\end{eqnarray}
leading to 
\begin{eqnarray}
m^3 \langle \overline{x} \ \overline{y} \ \overline{z} \rangle &= &m(m-1)(m-2) \langle x \rangle \langle y \rangle \langle z \rangle
\nonumber
 \\ \nonumber &+& m(m-1) \left(\langle xy \rangle \langle z \rangle +  \langle xz \rangle \langle y \rangle +  \langle yz \rangle \langle x \rangle \right) \\
 &+& m \langle xyz\rangle. \label{EqExpAv}
\end{eqnarray}
The above relations can expressed as one matrix equation
\begin{equation}\begin{pmatrix}
\langle \overline{xyz}\rangle\\
\langle \overline{xy} \ \overline{z} \rangle\\
\langle \overline{xz} \ \overline{y} \rangle\\
\langle \overline{yz} \ \overline{x} \rangle\\
\langle \overline{x} \ \overline{y} \ \overline{z} \rangle\\
\end{pmatrix}
=  \begin{pmatrix}
1 & 0 & 0 & 0 & 0 \\
\frac{m_1}{m^2} & \frac{m_2}{m^2} & 0 & 0 & 0 \\
\frac{m_1}{m^2} & 0 & \frac{m_2}{m^2} & 0 & 0 \\
\frac{m_1}{m^2} & 0 & 0 & \frac{m_2}{m^2} & 0 \\
\frac{m_1}{m^3}& \frac{m_2}{m^3} & \frac{m_2}{m^3} & \frac{m_2}{m^3} & \frac{m_3}{m^3}\\
\end{pmatrix} \cdot
\begin{pmatrix}
\langle xyz \rangle \\
\langle xy \rangle \langle z \rangle \\
\langle xz \rangle \langle y \rangle \\
\langle yz \rangle \langle x \rangle \\
\langle x \rangle \langle y \rangle \langle z \rangle\\
\end{pmatrix},
\end{equation}
where we defined $m_n = \prod_{i=1}^{n} (m-i+1)$. The $5\times5$ matrix will be denoted by $A_3$ in the following.  
The cumulant $C_3$ is then expressed as an statistical average of products of sample means
\begin{equation}
C_3(x,y,z)
=
\left\langle
\begin{pmatrix}
1 \\
-1 \\
-1 \\
-1 \\
2
\end{pmatrix}^T\cdot
A_3^{-1} \cdot
\begin{pmatrix}
 \overline{xyz} \\
 \overline{xy} \ \overline{z} \\
 \overline{xz} \ \overline{y} \\
 \overline{yz} \ \overline{x} \\
\overline{x} \ \overline{y} \ \overline{z}\\
\end{pmatrix}
\right\rangle \label{C3matrix}
\end{equation}
where the coefficients of the first vector follow from the RHS of Eq. (\ref{C_3}). 
The inversion of $A_3$ poses no problem and was performed by a computer algebra system.


It follows that the estimator $c_3(x,y,z)$ of the multivariate cumulant $C_3$ is given by the expression in the angle brackets which after evaluation yields
\begin{eqnarray}
c_3^{\rm (a)}(x,y,z) &=& \frac{m^2}{(m-1)(m-2)} \nonumber \\
\nonumber &&  \times( \overline{xyz} - \overline{xy} \ \overline{z} - \overline{xz} \ \overline{y} - \overline{yz} \ \overline{x} + 2 \overline{x} \ \overline{y} \ \overline{z} ) \nonumber \\
& = &  \frac{m^2}{(m-1)(m-2)} \overline{(x-\overline{x})(y-\overline{y})(z-\overline{z})}.
\end{eqnarray}
The estimator is defined for all sample sizes $m\geq3$ and will produce correct unbiased estimators in contrast to the 
'natural' estimator without the prefactor $\frac{m^2}{(m-1)(m-2)}$.
 The requirement $m \geq 3$ for the sample size is consistent with the fact that skewness cannot be determined from only two samples.\\
 Setting all variables equal to $x$
\begin{eqnarray}
		c_3^{\rm (b)}(x) & = & \frac{m^2}{(m-1)(m-2)} ( \overline{x^3} - 3\overline{x^2} \ \overline{x} + 2 \overline{x} \ \overline{x} \ \overline{x} ) \nonumber \\
		& = & \frac{m^2}{(m-1)(m-2)}  \overline{(x-\overline{x})^3} 		
\end{eqnarray} 
the third order k-statistics $k_3$ of Fisher is recoved \cite{FisherPLMS1928}. 

Several special cases that may be interesting for applications follow immediately from an adapted version of (\ref{C3matrix}). 
For  $\langle x \rangle = 0$ we find after replacing the first vector in (\ref{C3matrix}) by $(1,-1,-1,0,0)$
\begin{equation}
c_3^{\rm (e)}(x,y,z) =  \frac{1}{m-1} \left((m+1)\overline{xyz} - m\left(\overline{xy} \ \overline{z} + \overline{xz} \ \overline{y}\right)\right),
\end{equation}
which is valid for $m \ge 2$.  \\
For $\langle x \rangle = \langle y \rangle = 0$ we find
\begin{equation}
c_3^{\rm (f)}(x,y,z) =  \frac{m}{m-1} \left(\overline{xyz} - \overline{xy} \ \overline{z}\right).
\end{equation}
For all variables being average free ($\langle x \rangle = \langle y \rangle = \langle z \rangle = 0$), the result simplifies to
\begin{equation}
c_3^{\rm (d)}(x,y,z) =  \overline{xyz}.
\end{equation}
Another interesting case can be derived for two uncorrelated (not necessarily average-free) variables $x,y$ with $C_2(x,y)=0$. 
After replacing the first vector in Eq. (\ref{C3matrix}) by $(1,0,-1,-1,1)$ we find
\begin{equation}
	c_3^{\rm (i)}(x,y,z) = \frac{m}{m-2} \left(\overline{xyz} - \overline{xz} \ \overline{y} - \overline{yz} \ \overline{x} + \overline{x} \ \overline{y} \ \overline{z}\right).
\end{equation}
The estimator $c_3^{\rm (d)}$ simplifies to 
\begin{equation}
c_3^{\rm (c)}(x) =  \overline{x^3}
\end{equation}
for $x = y = z$.
A complex estimator for $x = a$, $y = a^*$, $z = b$ with $\langle a \rangle = 0$ may find application for the calculation of 
third order polyspectra (compare Section \ref{sec:c4c}).
We find
\begin{equation}
  c_3^{\rm (cd)}(a,b) = \frac{m}{m-1}(\overline{a a^* b} - \overline{a a^*}\, \overline{b}).
\end{equation}
\section{Estimators for $C_4$}
\label{sec:c4}
The derivation of unbiased estimators for $C_4$ follows the same scheme as above. We immediately find
\begin{equation}
\langle \overline{xyzw} \rangle = \langle xyzw \rangle. 
\end{equation}
The term  
\begin{equation}
	\langle \overline{xyz} \ \overline{w} \rangle = \frac{1}{m^2} \sum_{i,j}^{m} \langle x_i y_i z_i w_j \rangle
\end{equation}
exhibits the multiplicities 
\begin{eqnarray}
i=j  & \textrm{$m$-times} \nonumber  \\
 i\neq j   & \textrm{$m(m-1)$-times}. \nonumber
\end{eqnarray}
which yields
\begin{equation}
m^2 \langle \overline{xyz}  \ \overline{w} \rangle = \left[m(m-1) \langle xyz \rangle \langle w \rangle + m \langle xyzw \rangle \right].
\end{equation}
Corresponding expressions hold for $\overline{xyw} \ \overline{z}$, $\overline{xzw} \ \overline{y}$ , and $\overline{yzw} \ \overline{x}$.
The term 
\begin{equation}
\langle \overline{xy} \ \overline{zw} \rangle = \frac{1}{m^2} \sum_{i,j}^{m} \langle x_i y_i z_j w_j \rangle
\end{equation}
has the same multiplicities as $\langle \overline{xyz} \ \overline{w} \rangle$ and leads us to
\begin{equation}
m^2 \langle \overline{xy} \ \overline{zw} \rangle =  \left[m(m-1) \langle xy \rangle \langle zw \rangle + m \langle xyzw \rangle \right].
\end{equation}
Terms with structure $\overline{xy} \ \overline{z} \ \overline{w}$ can be written as 
\begin{equation}
\langle \overline{xy} \ \overline{z} \ \overline{w} \rangle = \frac{1}{m^3} \sum_{i,j,k}^{m} \langle x_i y_i z_j w_k \rangle
\end{equation}
with the multiplicities
\begin{eqnarray}
\nonumber i=j=k  \ &\textrm{$m$-times} \\
\nonumber i=j\neq k  \ &\textrm{$m(m-1)$-times} \\
\nonumber i\neq j = k \ &\textrm{$m(m-1)$-times} \\
\nonumber i = k \neq j \ &\textrm{$m(m-1)$-times} \\
\nonumber i \neq j\neq k  \ &\textrm{ $m(m-1)(m-2)$-times } 
\end{eqnarray}
resulting in
\begin{eqnarray}
m^3\langle \overline{xy} \ \overline{z} \ \overline{w} \rangle =  &&m(m-1)(m-2) \langle xy \rangle \langle z \rangle \langle w \rangle \\
\nonumber &+& m(m-1)( \langle xyz \rangle \langle w \rangle + \langle xyw \rangle \langle z \rangle ) \\
\nonumber &+& m(m-1) \langle xy \rangle \langle zw \rangle + m \langle xyzw \rangle.
\end{eqnarray}
Corresponding expressions hold for $ \langle\overline{xz} \ \overline{y} \ \overline{w}\rangle$, $\langle\overline{xw} \ \overline{y} \ \overline{z} \rangle$, $\langle\overline{yz} \ \overline{x} \ \overline{w}\rangle$, $\langle\overline{yw} \ \overline{x} \ \overline{z}\rangle$ and $\langle\overline{zw} \ \overline{x} \ \overline{y}\rangle$. Finally, the term
\begin{equation}
\langle \overline{x} \ \overline{y} \ \overline{z} \ \overline{w} \rangle= \frac{1}{m^4} \sum_{i,j,k,l}^{m} \langle x_i y_j z_k w_l \rangle
\end{equation}
exhibits the following multiplicities 
\begin{eqnarray}
\nonumber i=j=k=l  \ &\textrm{$m$-times} \\
\nonumber i=j=k\neq l  \ &\textrm{$m(m-1)$-times with 4 realizations} \\
\nonumber i=j\neq k=l  \ &\textrm{$m(m-1)$-times with 3 realizations} \\
\nonumber i=j \neq k\neq j  \ &\textrm{ $m(m-1)(m-2)$-times with 6 realizations} \\
\nonumber i\neq j\neq k \neq l \ &\textrm{$m(m-1)(m-2)(m-3)$-times}, 
\end{eqnarray}
where e.g. '4 realizations'  means explicitly $i = j = k \neq l$, $i = j = l \neq k$, $j = k = l \neq i$, and $i = k = l \neq j$.
This leads to
\begin{eqnarray}
 \nonumber m^4 \langle \overline{x} \ \overline{y} \ \overline{z} \ \overline{w} \rangle = &&m(m-1)(m-2)(m-3) \langle x \rangle \langle y \rangle \langle z \rangle \langle w \rangle  \\
\nonumber &+& m(m-1)(m-2) (\langle xy \rangle \langle z \rangle \langle w \rangle + \textrm{5 o.p.}) \\
\nonumber &+& m(m-1) (\langle xyz \rangle \langle w \rangle  +  \textrm{3 o.p.} ) \\
\nonumber &+& m(m-1) (\langle xy \rangle \langle zw \rangle +  \textrm{2 o.p.} ) \\
 &+& m \langle xyzw\rangle, \label{eq:c_example}
\end{eqnarray}
where 'o.p.' means other permutations of the variables in e.g. $\langle xy \rangle \langle z \rangle \langle w \rangle$ that give rise to
(non-identical) terms like  $\langle xz \rangle \langle y \rangle \langle w \rangle$. 

The relation of means and products of expected statistical averages can be written as a matrix equation
\begin{equation}
 \langle \vec{\mu} \rangle = A_4 \vec{p}
\end{equation}
where
\begin{equation}
\vec{\mu} = 
\begin{pmatrix}
 \overline{xyzw}\\
 \overline{xyz}\ \overline{w}\\
 \overline{xyw}\ \overline{z}\\
 \overline{xzw}\ \overline{y}\\
 \overline{yzw}\ \overline{x}\\
 \overline{xy} \ \overline{zw}\\
 \overline{xz} \ \overline{yw}\\
 \overline{xw} \ \overline{yz}\\
 \overline{xy} \ \overline{z} \ \overline{w} \\
 \overline{xz} \ \overline{y} \ \overline{w} \\
 \overline{xw} \ \overline{y} \ \overline{z} \\
 \overline{yz} \ \overline{x} \ \overline{w} \\
 \overline{yw} \ \overline{x} \ \overline{z} \\
 \overline{zw} \ \overline{x} \ \overline{y} \\
 \overline{x} \ \overline{y} \ \overline{z} \ \overline{w} 
\end{pmatrix}; \quad\quad
\vec{p} = 
\begin{pmatrix}
\langle xyzw \rangle \\
\langle xyz \rangle \langle w \rangle \\
\langle xyw \rangle \langle z \rangle \\
\langle xzw \rangle \langle y \rangle \\
\langle yzw \rangle \langle x \rangle \\
\langle xy \rangle \langle zw \rangle \\
\langle xz \rangle \langle yw\rangle \\
\langle xw \rangle \langle yz \rangle \\
\langle xy \rangle\langle z \rangle\langle w \rangle \\
\langle xz \rangle\langle y \rangle\langle w \rangle \\
\langle xw \rangle\langle y \rangle\langle z \rangle \\
\langle yz \rangle\langle x \rangle\langle w \rangle \\
\langle yw \rangle\langle x \rangle\langle z \rangle \\
\langle zw \rangle\langle x \rangle\langle y \rangle \\
\langle x \rangle \langle y \rangle \langle z \rangle \langle w \rangle \\
\end{pmatrix}
\end{equation}
with the coefficient matrix $A_4$ given in Table \ref{A4}. 
\begin{table*}
	\normalsize
	\caption{\label{A4}  Matrix $A_4$}
	\renewcommand*{\arraystretch}{1.2}
	\begin{equation}
	A_4 = 
	\begin{pmatrix}
	1 & 0 & 0 & 0 & 0 & 0 & 0 & 0 & 0 & 0 & 0 & 0 & 0 & 0 & 0\\
	\frac{m_1}{m^2} & \frac{m_2}{m^2} & 0 & 0 & 0 & 0 & 0 & 0 & 0 & 0 & 0 & 0 & 0 & 0 & 0\\
	\frac{m_1}{m^2} & 0 & \frac{m_2}{m^2} & 0 & 0 & 0 & 0 & 0 & 0 & 0 & 0 & 0 & 0 & 0 & 0\\
	\frac{m_1}{m^2} & 0 & 0 & \frac{m_2}{m^2} & 0 & 0 & 0 & 0 & 0 & 0 & 0 & 0 & 0 & 0 & 0\\
	\frac{m_1}{m^2} & 0 & 0 & 0 & \frac{m_2}{m^2} & 0 & 0 & 0 & 0 & 0 & 0 & 0 & 0 & 0 & 0\\
	\frac{m_1}{m^2} & 0 & 0 & 0 & 0 & \frac{m_2}{m^2} & 0 & 0 & 0 & 0 & 0 & 0 & 0 & 0 & 0\\
	\frac{m_1}{m^2} & 0 & 0 & 0 & 0 & 0 & \frac{m_2}{m^2} & 0 & 0 & 0 & 0 & 0 & 0 & 0 & 0\\
	\frac{m_1}{m^2} & 0 & 0 & 0 & 0 & 0 & 0 & \frac{m_2}{m^2} & 0 & 0 & 0 & 0 & 0 & 0 & 0\\
	\frac{m_1}{m^3} & \frac{m_2}{m^3} & \frac{m_2}{m^3} & 0 & 0 & \frac{m_2}{m^3} & 0 & 0 & \frac{m_3}{m^3}& 0 & 0 & 0 & 0 & 0 & 0\\
	\frac{m_1}{m^3} & \frac{m_2}{m^3} & 0 & \frac{m_2}{m^3} & 0 & 0 & \frac{m_2}{m^3} & 0 & 0 & \frac{m_3}{m^3}& 0 & 0 & 0 & 0 & 0\\
	\frac{m_1}{m^3} & 0 & \frac{m_2}{m^3} & \frac{m_2}{m^3} & 0 & 0 & 0 & \frac{m_2}{m^3} & 0 & 0 & \frac{m_3}{m^3}& 0 & 0 & 0 & 0\\
	\frac{m_1}{m^3} & \frac{m_2}{m^3} & 0 & 0 &  \frac{m_2}{m^3} & 0 & 0 & \frac{m_2}{m^3} & 0 & 0 & 0 &  \frac{m_3}{m^3}& 0 & 0 & 0 \\
	\frac{m_1}{m^3} & 0 & \frac{m_2}{m^3} & 0 & \frac{m_2}{m^3} & 0 & \frac{m_2}{m^3} & 0 & 0 & 0 & 0 & 0 &\frac{m_3}{m^3}& 0 & 0\\
	\frac{m_1}{m^3} & 0 & 0 &\frac{m_2}{m^3} & \frac{m_2}{m^3} & \frac{m_2}{m^3} & 0 & 0 & 0 & 0 & 0 & 0 & 0& \frac{m_3}{m^3}& 0 \\
	\frac{m_1}{m^4}& \frac{m_2}{m^4} & \frac{m_2}{m^4} & \frac{m_2}{m^4} & \frac{m_2}{m^4} & \frac{m_2}{m^4} & \frac{m_2}{m^4} & \frac{m_2}{m^4} & \frac{m_3}{m^4} & \frac{m_3}{m^4} & \frac{m_3}{m^4} & \frac{m_3}{m^4} & \frac{m_3}{m^4} & \frac{m_3}{m^4} & \frac{m_4}{m^4}\\
	\end{pmatrix} \nonumber
	\end{equation}
	\hrulefill
		\label{A4}
\end{table*}
The cumulant $C_4$ [Eq.  (\ref{C_4})] can now be expressed as 
\begin{equation}
C_4(x,y,z,w)
= \langle \vec{\gamma}A_4^{-1}\vec{\mu} \rangle
\end{equation}
where $\vec{\gamma}^{\rm T} = (1,-1,-1,-1,-1,-1,-1,2,2,2,2,2,2,-6)$.
The full unbiased estimator for the fourth-order cumulant is found in the angle brackets as
\begin{eqnarray}
 c_4^{\rm (a)}(x,y,z,w) =&& \frac{m^2}{(m-1)(m-2)(m-3)} \times \\
&& \nonumber \Big[(m+1) \overline{xyzw} - (m+1) \left(\overline{xyz} \ \overline{w} + \textrm{ 3 o.p.}\right) \Big. \\ 
&&\nonumber \Big. -  (m-1) \left(\overline{xy} \ \overline{zw} + \textrm{2 o.p.} \right)  \Big. \\
&&\nonumber \Big.+ 2 m \left(\overline{xy} \ \overline{z} \ \overline{w} +\textrm{5 o.p.}\right) - 6m \overline{x} \ \overline{y} \ \overline{z} \ \overline{w} \Big.]
\end{eqnarray}
or more compactly written as
\begin{eqnarray}
 c_4^{\rm (a)}(x,y,z,w) =&& \frac{m^2}{(m-1)(m-2)(m-3)}  \\
&& \nonumber \times\Big[(m+1) \overline{(x-\overline{x})(y-\overline{y})(z-\overline{z})(w-\overline{w})} \Big. \\
&& \nonumber \Big. -  (m-1) \left(\overline{(x-\overline{x})(y-\overline{y})} \right. \\ 
&&   \nonumber \left. \times \overline{(z-\overline{z})(w -\overline{w})}  + \textrm{2 o.p.} \right) \Big].
\end{eqnarray}
Unlike in the case of $c_2$ and $c_3$, there no longer is a single common prefactor for all terms. 
Several special cases for $c_4$ will be discussed in the following. In the case of all variables being equal to $x$
we find
  \begin{eqnarray}
c_4^{\rm (b)}(x)  &= &\frac{m^2}{(m-1)(m-2)(m-3)} \times \\
 \nonumber && \left[(m+1) \overline{x^4} - 4 (m+1) \overline{x^3} \ \overline{x} \right. \\
 \nonumber && \left.- 3 (m-1) \overline{x^2} \ \overline{x^2} + 12 m \overline{x^2} \ \overline{x} \ \overline{x} - 6m \overline{x}^4\right] 
 \nonumber \\
  & = & \frac{m^2[(m+1) \overline{(x-\overline{x})^4} - 3(m-1)\overline{(x-\overline{x})^2}^2]}{(m-1)(m-2)(m-3)}  \nonumber \\
  &&
 \end{eqnarray}
which agrees with the fourth order k-statistics $k_4$ \cite{FisherPLMS1928}.
For four variables and $\langle x \rangle = 0$ we find
\begin{eqnarray}
c_4^{\rm (e)}(x,y,z,w)=&& \frac{1}{\left(m-1\right)\left(m-2\right)} \times \\ 
\nonumber &&\left[\left(m+1\right)\left(m+2\right) \overline{xyzw} \right.\\
\nonumber &&\left.-m \left(\left(m+2\right) \left(\overline{xyz} \ \overline{w} + \overline{xyw} \ \overline{z} + \overline{xwz} \ \overline{y}\right) \right. \right.\\
\nonumber && \left. \left. + m\left(\overline{xy} \ \overline{zw} + \overline{xz} \ \overline{yw} + \overline{xw} \ \overline{yz} \right) \right. \right. \\
\nonumber && \left. \left. - 2m \left(\overline{xy} \ \overline{z} \ \overline{w} + \overline{xz} \ \overline{y} \ \overline{w} + \overline{xw} \ \overline{y} \ \overline{z}\right)\right) \right].
\end{eqnarray}
For $\langle x \rangle = \langle y \rangle = 0$ we find
\begin{eqnarray}
c_4^{\rm (f)}(x,y,z,w) = &&\frac{1}{\left(m-1\right)\left(m-2\right)} \times \\
\nonumber && \left(\left(m^2+2m-4\right)\overline{xyzw} \right.\\
\nonumber &&\left. -m \left(m \left(\overline{xyz} \ \overline{w} + \overline{xyw} \ \overline{z}\right) + m \overline{xy} \ \overline{zw}\right. \right. \\
\nonumber &&\left. \left. + \left(m-2\right) \left( \overline{xz} \ \overline{yw} + \overline{xw} \ \overline{yz} \right) \right. \right. \\
\nonumber &&\left. \left.- 2m \overline{xy} \ \overline{z} \ \overline{w}\right) \right).
\end{eqnarray}
For $\langle x \rangle = \langle y \rangle = \langle z \rangle = 0$ we find
\begin{eqnarray}
&&c_4^{\rm (g)}(x,y,z,w) = \frac{1}{m-1} \times \\ 
\nonumber && \left(\left(m+3\right) \overline{xyzw} - m \left(\overline{xyz} \ \overline{w} + \overline{xy} \ \overline{zw} + \overline{xz} \ \overline{yw} + \overline{xw} \ \overline{yz} \right)\right).
\end{eqnarray}
And for all variables being average-free we obtain
\begin{eqnarray}
\label{C4zeromean}
&&c_4^{\rm (d)}(x,y,z,w) = \frac{1}{m-1} \times \\
\nonumber &&\left(\left(m+2\right) \overline{xyzw} - m \left( \overline{xy} \ \overline{zw} + \overline{xz} \ \overline{yw} + \overline{xw} \ \overline{yz} \right)\right).
\end{eqnarray}
The problem of finding $c_4^{\rm (d)}(x,y,z,w)$ was stated in \cite{BlagouchineIEEE2009} but left unsolved.
Also for average-free variables there is no common prefactor for all terms. For one variable $x$ with zero mean we recover the 
known result \cite{mansourCONF1998,BlagouchineIEEE2009}
\begin{equation}
c_4^{\rm (c)}(x) = \frac{1}{m-1} \left(\left(m+2\right) \overline{x^4} - 3m \overline{x^2} \ \overline{x^2})\right).
\end{equation} 
For four average free variables with vanishing pairwise second order cumulant $C_2(x,y) = 0$,  $C_2(x,z) = 0$, etc. we obtain the simple
estimator
\begin{equation}
c_4^{\rm (j)}(x) =  \overline{x y z w}.
\end{equation} 


\section{Estimators of $C_4$ for complex variables}
\label{sec:c4c}
For complex variables $a$, $b$, $c$, $d$ the above estimators can be used analogously. 
A few additional special cases are however interesting for calculating fourth order polyspectra.
 In Ref.  \cite{starosielecRSI2010} the authors obtained a fourth order spectrum from estimating a cumulant
 $C_4(a_\omega,a^*_{\omega},a_{\omega'},a^*_{\omega'})$. The variables $a_{\omega}$ were Fourier coefficients obtained from Fast Fourier Transformations of a stochastic signal.
The variables $a_\omega$ and $a_{\omega'}$ always exhibit a random complex phase except for $\omega = 0$ where a
constant offset may appear if the initial stochastic signal is not average-free.
 The following special cases are therefore highly relevant for $C_4$ estimations.
The cumulant $C_4(a,a^*,b,b^*)$ is for $\langle a\rangle = \langle b \rangle = 0$ and $C_2(a,b) = C_2(a,b^*) = 0$ given by
\begin{equation}
 c_4^{\rm (ca)}(a,b) = \frac{m}{m-1}(\overline{aa^*bb^*} - \overline{aa^*} \ \overline{bb^*}). 
\end{equation}
The estimator $c_4^{\rm (ca)}(a,b)$ had previously been derived by Starosielec  \cite{starosielecBOOK2012}
and found application in \cite{starosielecRSI2010}.

If $C_2(a,b) = 0$ can be assumed but $C_2(a,b^*) = 0$  cannot (e.g. if $a = b$), the estimator
\begin{eqnarray}
 c_4^{\rm (cb)}(a,b) & = & \frac{1}{m-1} \left(\left(m+1\right) \overline{aa^*bb^*}\right. \nonumber \\
 & &\left. - m \left( \overline{a a^*} \ \overline{bb^*}  + \overline{ab^*} \ \overline{a^*b} \right)\right)
\end{eqnarray}
 has to be used. 
 
The authors of  \cite{starosielecRSI2010} used $c_4^{\rm (ca)}(a_\omega,a_{\omega'})$ with $m = 2$ to estimate $C_4$ for  all frequency pairs $\omega$, $\omega'$. While this is correct for cases where $\omega \neq \omega'$ the estimator
  $c_4^{\rm (cb)}(a_\omega,a_{\omega'})$ should have been used for $\omega = \omega'$
  since $C_2(a_\omega,a_{\omega'}^*) \neq 0$ in that case. In case of a purely Gaussian signal the authors
  found falsely a strong contribution for  $\omega = \omega'$ while  $C_4(a_\omega,a^*_{\omega},a_{\omega},a^*_{\omega})$ 
  is in fact strictly zero in such a case.
  
  For $\langle a \rangle \neq 0$, $\langle b \rangle = 0$, $C_2(a,b) = 0$, and $C_2(a,b^*) = 0$ we obtain
\begin{eqnarray}
 c_4^{\rm (cc)}(a,b) & = & \frac{m^2}{(m-1)(m-2)} \left( \overline{aa^*bb^*}  -  \overline{a^* b b^*}\,\overline{a}\right. \nonumber \\
 & & \left. - \overline{a b b^*} \, \overline{a^*}
 - \overline{a a^*} \,\overline{b b^*} + 2 \overline{b b^*} \, \overline{a}\,\overline{a^*} \right).
\end{eqnarray} 
This estimator $c_4^{\rm (cc)}(a,b)$ may find application for a polyspectrum if $a = a_{\omega}$ with $\omega \approx 0$.
If $C_4(a_{\omega_1},a_{\omega_2},a_{\omega_3},a_{\omega_4})$ needs to be estimated with all frequency pairs $\omega_j + \omega_i \neq 0$, the estimator $c_4^{\rm (j)}$ should be used. 
\section{The Variances of Estimators}
\label{sec:var}
\begin{figure}
	\centering
	\includegraphics[width=7cm]{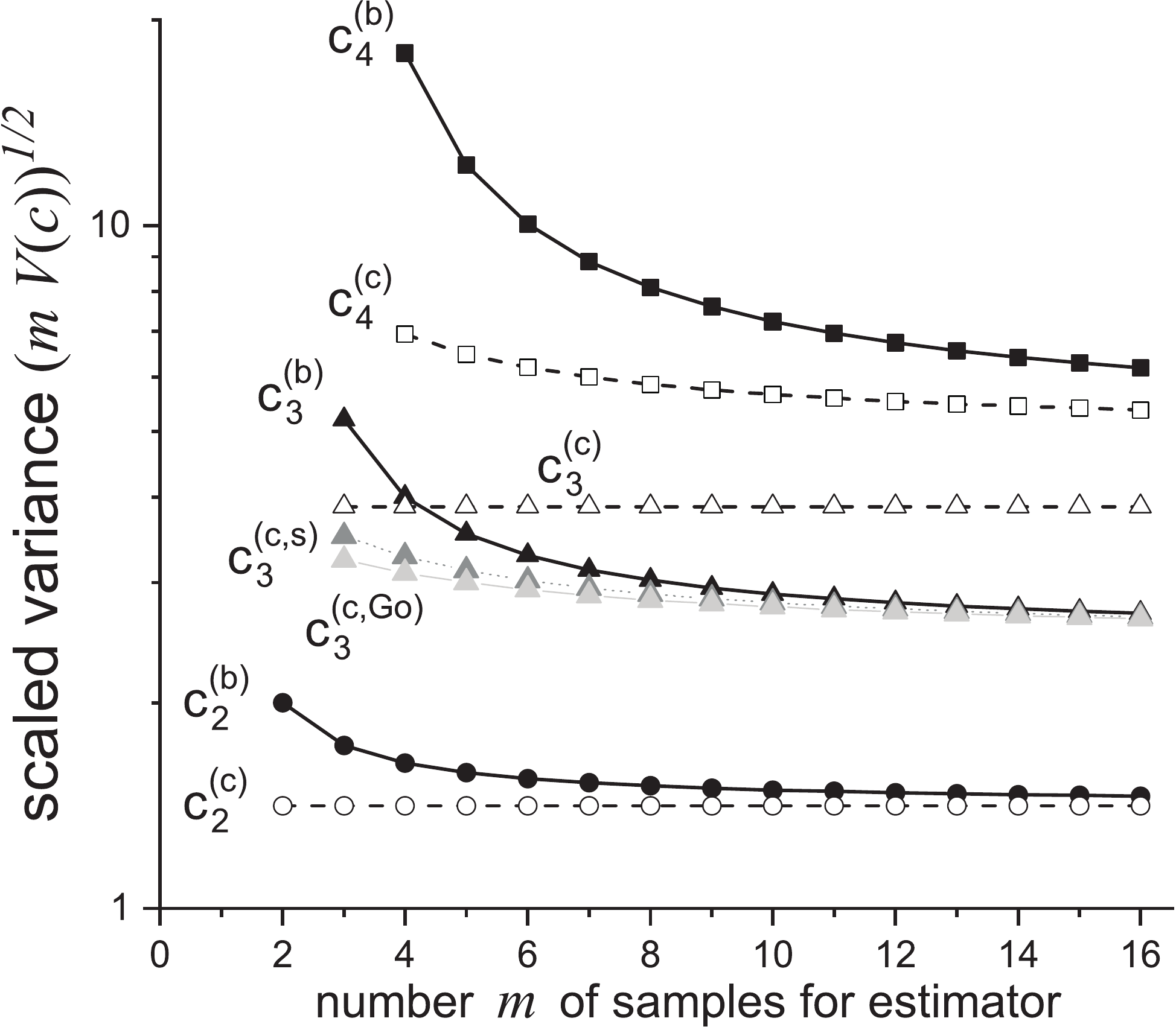}
	\caption{The performance of different estimators is revealed by the scaled variance $\sqrt{m V(c)}$. The variance is calculated for 
	Gaussian processes with $\langle x \rangle = 0$ and  $\langle x^2 \rangle = 1$.  The reduced estimators $c_2^{\rm (c)}$
	and $c_4^{\rm (c)}$  (broken lines) are superior over the full estimators $c_2^{\rm (b)}$
	and $c_4^{\rm (b)}$ (solid lines). Surprisingly, the full estimator $c_3^{\rm (b)}$ performs better than its reduced counterpart $c_3^{\rm (c)}$ for $m\geq 5$. The ultimately best estimator for $C_3(x,x,x)$ is $c_3^{\rm (c,Go)}$.}
	\label{fig:EstimatorPerformance}
\end{figure}
In a real-world application an approximation $\tilde{C}$ of a cumulant $C$ is calculated from a limited number $M$ of estimates $c$ via
\begin{equation}
 \tilde{C} = \langle c \rangle_{M}.
\end{equation} 
The variance $\sigma^2$ of $\tilde{C}$
\begin{equation}
  \sigma^2 = \frac{\langle c^2 \rangle - \langle c \rangle^2}{M} = \frac{m(\langle c^2 \rangle - \langle c \rangle^2)}{m M}. 
\end{equation} 
is a measure of how accurate $C$ can be determined from $m M$ samples of the random variables.
Regarding this, the scaled variance
\begin{equation}
   m V(c) =m(\langle c^2 \rangle - \langle c \rangle^2) 
\end{equation} 
is a sensible measure for comparing the performance of different estimators $c$ with varying $m$.
We calculate $V(c_n)$ for a selection of estimators in Appendix \ref{app:var}. The results can be expressed in terms of cumulants 
of up to order $2n$. In general $V(c_n)$ can like $C_n$ only be estimated from samples of the random process.
If the process is however dominated by Gaussian noise, all cumulants of order three or higher no longer contribute to $V(c_n)$ and 
the expressions for $V(c_n)$ greatly simplify.

In the following, we compare univariate estimators where $x$  is dominated by a Gaussian  
contribution and $\langle x \rangle = 0$. Such processes often appear in physics where the use of an AC-coupled amplifier leads to average free time-series. 
The arrival of photons from a laser in a detector is known to be Poisson-distributed. The signal is amplified and can be sampled to yield a series of data points that should exhibit an almost Gaussian, but still slightly asymmetric distribution of values centered around zero. Figure \ref{fig:EstimatorPerformance} shows the square root of the scaled variance, i.e. $\sqrt{m V(c)}$, for six different estimators for a random Gaussian process $x$ with $\langle x \rangle = 0$ and   $\langle x^2 \rangle = 1$.
The value of  $\sqrt{m V(c)}$ for the full estimators $c_2^{\rm (b)}$, $c_3^{\rm (b)}$, and $c_4^{\rm (b)}$, 
which do not require $\langle x \rangle =0$, are plotted as solid lines for increasing $m$. The reduced estimators $c_2^{\rm (c)}$, $c_3^{\rm (c)}$, and $c_4^{\rm (c)}$, which require $\langle x \rangle = 0$, are shown as dashed lines. The full expressions of $V(c)$ for all estimators can be found in Appendix \ref{app:var}. 
In general, the estimators exhibit a larger $\sqrt{m V(c)}$ for increasing order. A factor of 10 between the fourth and second order estimator results in practice in 100 times more samples that are required in the $C_4$ case to obtain a similar noise level of the estimate as compared to the $C_2$ case. 

The reduced estimators perform always better for $C_2$ and $C_4$. Surprisingly, we find that  the reduced estimator $c_3^{\rm (c)}$ for $C_3$ performs worse than the full estimator $c_3^{\rm (b)}$ for $m \geq 5$. Despite the knowledge of $x$ being average free, the reduced estimator yields no benefit and should in fact be avoided for practical purposes. 
\section{Gauss-optimal Estimators}
\label{sec:go}
Next, we show that even better estimators for $C_3(x,x,x)$ than $c_3^{\rm (b)}$ can be found.
The superposition $c_3^{\rm (g,s)} = \alpha c_3^{\rm (b)} + (1- \alpha) c_3^{\rm (c)}$ is an unbiased estimator of $C_3(x,x,x)$ for  $\langle x \rangle = 0$ . In case of a Gaussian process we find the scaled variance
\begin{equation}
  m V(c_3^{\rm (c,s)}) = \left(15 - 18 \alpha + \frac{9 m^2 - 9m + 6}{(m-1)(m-2)}\alpha^2\right)\langle x^2 \rangle^3
 \end{equation}
which assumes a minimal value of $6(3 m^2 + 6 m -4) \langle x^2 \rangle^3 /(3 m^2 + 6m -4)$ for $\alpha = (3 m^2 -9m +6)/(3 m^2 - 3 m +2)$. The minimal values for $\sqrt{  m V(c_3^{\rm (c,s)}) }$ for increasing $m$ are plotted in Figure \ref{fig:EstimatorPerformance} (dotted line) and show a significant improvement over the estimators $c_3^{(b)}$ and $c_3^{(c)}$. 
 
The utmost best estimator ist found by considering the most general estimator for $C_3(x,x,x)$ with $\langle x \rangle = 0$ 
\begin{equation}
c_3^{\rm (c,gen)} = 
\begin{pmatrix}
1\\
\alpha_1 \\
\alpha_2
\end{pmatrix}^{\rm T}
\begin{pmatrix}
1 & 0 & 0  \\
\frac{m_1}{m^2} & \frac{m_2}{m^2} & 0  \\
\frac{m_1}{m^3}& 3 \frac{m_2}{m^3} & \frac{m_3}{m^3}\\
\end{pmatrix}^{-1}
\begin{pmatrix}
 \overline{x^3}\\
\overline{x^2} \ \overline{x} \\
 \overline{x}^3 
\end{pmatrix}.
\end{equation}
The equation is very similar to (\ref{C3matrix}) where the matrix was adapted for the case of a single variable $x$.
We find $\langle  c_3^{\rm (c,gen)} \rangle = \langle x^3 \rangle + \alpha_1 \langle x^2 \rangle \langle x \rangle + \alpha_2 \langle x \rangle^3$ which means that for any non-zero $\alpha_1$ or $\alpha_2$ a zero is effectively added to $C_3$. While the expectation value of the estimator $\langle  c_3^{\rm (c,gen)} \rangle$ does not depend on $\alpha_i$, an optimized variance $V( c_3^{\rm (c,gen)})$ may be found for non-zero $\alpha_i$.
Here the estimator is optimized for a  Gaussian process $x$.
A fully analytic solution is possible since  $V(c_3^{\rm (c,gen)})$  is only quadratic in $\alpha_1$ and $\alpha_2$. The Gauss-optimal (Go) estimator
\begin{equation}
c_3^{\rm (c, Go)} = \overline{x^3} - \frac{3(m-1)}{m+1}\overline{x^2} \ \overline{x} \label{eq:c3c_Go}
\end{equation}
with 
\begin{equation}
m V(c_3^{\rm (c,Go)}) = \frac{6(m+4)\langle x^2 \rangle^3}{m+1}
\end{equation}
follows with the help of computer algebra. 
The optimal values of $\sqrt{m V(c_3^{\rm (c,Go)}) }$ for increasing $m$ are plotted in Figure \ref{fig:EstimatorPerformance} (light grey line). The estimator $c_3^{\rm (c,Go)}$ surpasses the performance of the three other estimators $c_3^{\rm (b)}$, $c_3^{\rm (c)}$, and $c_3^{\rm (c,s)}$.
 We are not aware that $c_3^{\rm (c,Go)}$ has been discovered before in literature.
A corresponding calculation for $C_2$ yields no improvement over $c_2^{\rm (c)}(x)$.
Similarly, a longer calculation shows for $C_4$ that the Gauss-optimal estimator is given by the reduced estimator $c_4^{(c)}(x)$ discussed 
previously by Blagouchine and Moreau \cite{BlagouchineIEEE2009}.

\begin{figure}
	\centering
	\includegraphics[width=7cm]{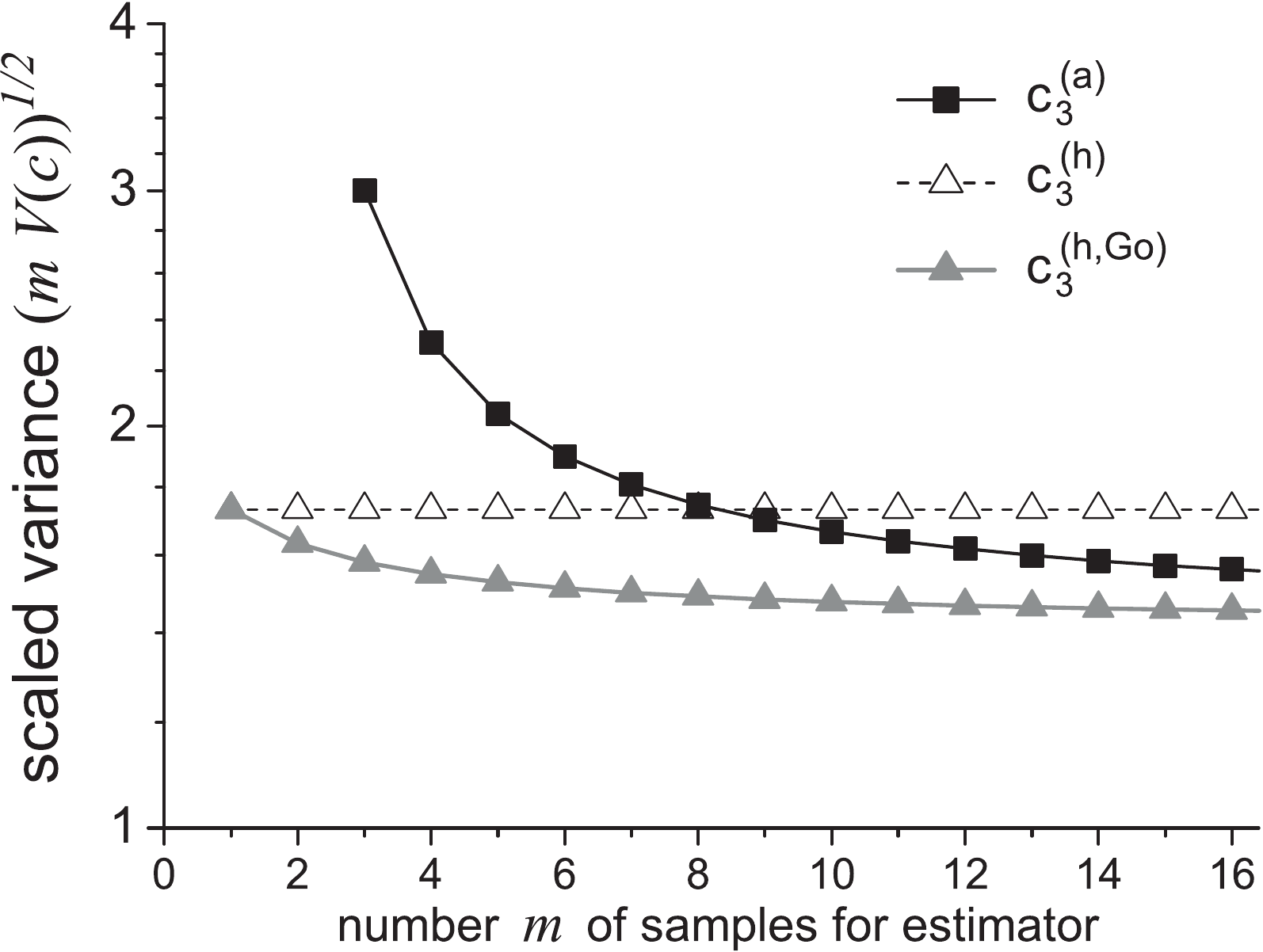}
	\caption{Comparison of the performance $\sqrt{m V(c)}$ of three different estimators for  $C_3(x,x,z)$ in the case of two independent average-free Gaussian processes $x$ and $z$ with $\langle x^2 \rangle = \langle z^2 \rangle = 1$. The new Gauss-optimal estimator $c_3^{\rm (h,Go)}(x,z)$ performs for any $m$ samples better than both the reduced estimator $c_3^{\rm (h)}(x,z)$ and the full estimator $c_3^{\rm (a)}(x,x,z)$.}
	\label{fig:xxzEstimator}
\end{figure}

Next, we derive Gauss-optimal estimators for cumulants of more than one variable 
using the procedure above. The Gauss-optimal second order estimator for $C_2(x,y)$, where $\langle x \rangle = 0$ 
and $\langle y \rangle \neq 0$, is identical with $c_2^{\rm (d)}(x,y) = \overline{xy}$. The same holds true for $c_2^{\rm (e)}(x,y) = \overline{xy}$ when $\langle x \rangle = \langle y \rangle =  0$.

The reduced estimator for $C_3(x,x,z)$  for  $\langle x \rangle = \langle z \rangle = 0$ 
is 
\begin{equation}
 c_3^{\rm (h)}(x,z)  =  \overline{x^2 z}. 
\end{equation}
with
\begin{equation}  
    m V(c_3^{\rm (h)}(x,z)) = 3\langle x^2 \rangle^2 \langle z^2 \rangle.
\end{equation}
The Gauss-optimal version of that estimator 
\begin{equation}
  c_3^{\rm (h,Go)}(x,z)  = \frac{m+2}{m+1} \overline{x^2 z} - \frac{m}{m+1}\overline{x^2}\,\overline{z} 
  \end{equation}
exhibits an improved variance 
\begin{equation}  
    m V(c_3^{\rm (h,Go)}(x,z)) =  \frac{2(m+2)\langle x^2 \rangle^2 \langle z^2 \rangle}{m+1}.  
\end{equation}
which outperforms the reduced estimator $\overline{x^2 z}$ by a factor of up to 3/2 (compare also Figure \ref{fig:xxzEstimator}). 
This authors of \cite{gabelliNJP2013} investigated the current statistics of a quantum electronics device using 
estimator $c_3^{\rm (h)}(x,z)$ implemented in hard wired analog electronics. Similar experiments may in the future benefit from the use of 
 the Gauss-optimal estimator $c_3^{\rm (h,Go)}(x,z)$.
The general estimator $c_3^{\rm (a)}(x,x,z)$ performs for large $m$ almost as good as $c_3^{\rm (h,Go)}(x,z)$:
\begin{eqnarray}
 m V(c_3^{\rm (a)}(x,x,z))& = & \frac{2 m^2 \langle x^2 \rangle^2 \langle z^2 \rangle}{(m-1)(m-2)}. 
\end{eqnarray}
The estimator for $C_3(x,y,z)$
\begin{equation}
  c_3^{\rm (d)}(x,y,z)  = \overline{xyz}  
\end{equation}
is Gauss-optimal 
with
\begin{equation}  
  m V(c_3^{\rm (d)}(x,y,z)) = \langle x^2 \rangle \langle y^2 \rangle\langle z^2 \rangle.
\end{equation}
For comparison, the general estimator $c_3^{(a)}$ performs worse with
\begin{eqnarray}
m V(c_3^{\rm (a)}(x,y,z)) & = & \frac{m^2 \langle x^2 \rangle \langle y^2 \rangle\langle z^2 \rangle}{(m-1)(m-2)}.
\end{eqnarray}
Similarly,
we find that $c_3^{\rm (e)}(x,y,z)$, $c_3^{\rm (f)}(x,y,z)$ and $c_3^{\rm (cb)}(a,b)$ are Gauss-optimal with
\begin{eqnarray}
 m V(c_3^{\rm (e)}(x,y,z)) & = &   \frac{m+1}{m-1}\langle x^2 \rangle C_2(y,y) C_2(z,z) \\
 m V(c_3^{\rm (f)}(x,y,z)) & = &  \frac{ m \langle x^2 \rangle \langle y^2 \rangle C_2(z,z)}{m -1} \\
 m V(c_3^{\rm (cd)}(a,b) ) & = & \frac{2 m\langle a a^* \rangle^2 C_2(b,b^*)}{m -1}.
\end{eqnarray}  

The fourth order case for several variables is extremely intricate and will not be discussed here. We suspect, however, that 
the reduced fourth order estimators are always Gauss-optimal like the univariate case $c_4^{(c)}(x)$ (see above).

\section{Numerical example for estimating $C_3(x,x,x)$}
\label{sec:num}
\begin{figure}
	\centering
	\includegraphics[width=7cm]{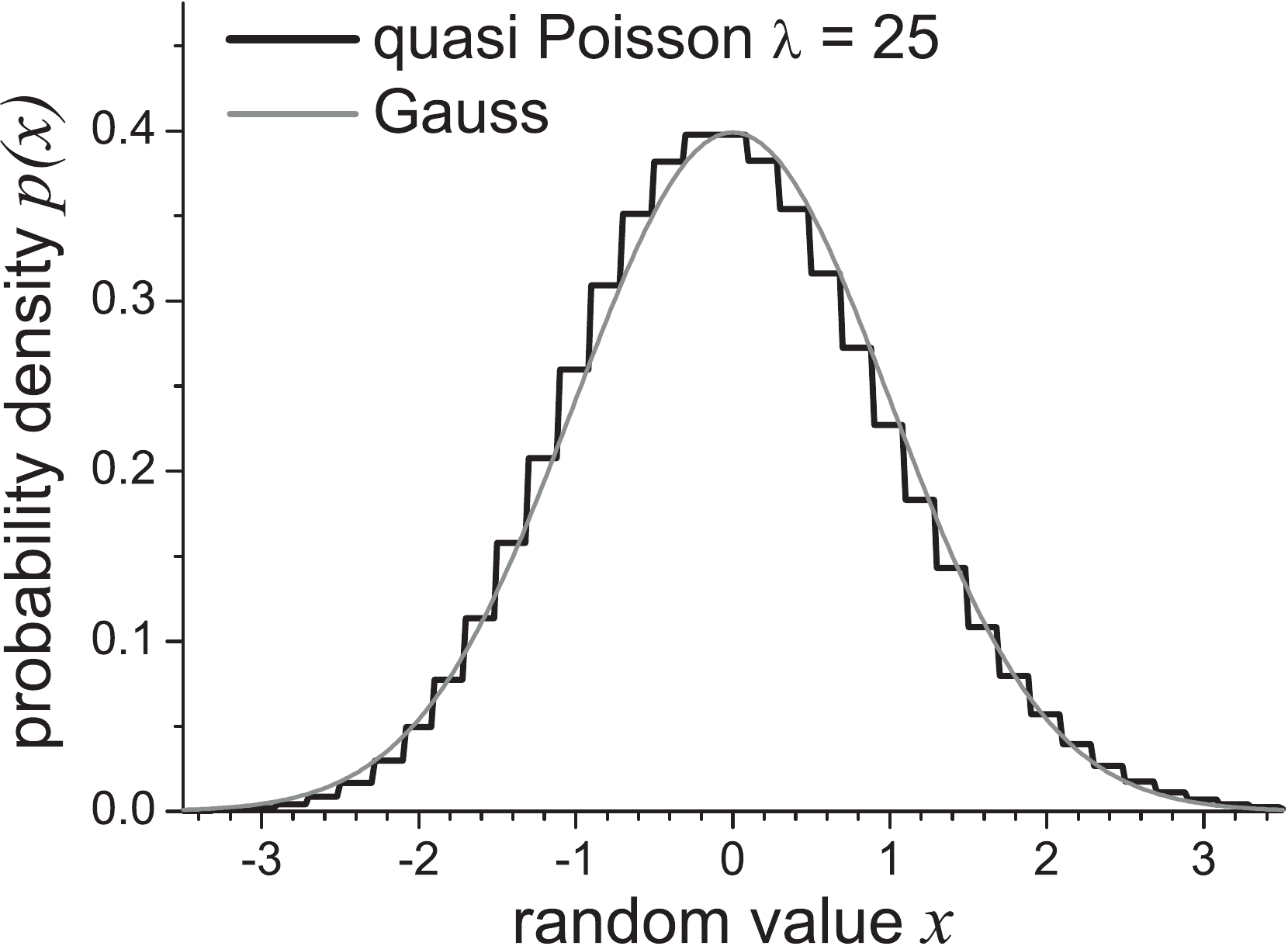}
	\caption{Comparison of an average-free quasi Poisson distribution with $\lambda = 25$ and a Gauss distribution.
The slight asymmetry of the quasi Poisson distribution causes a non-zero $C_3(x,x,x)$.}
	\label{fig:QuasiPoisson}
\end{figure}
\begin{figure}
	\centering
	\includegraphics[width=7cm]{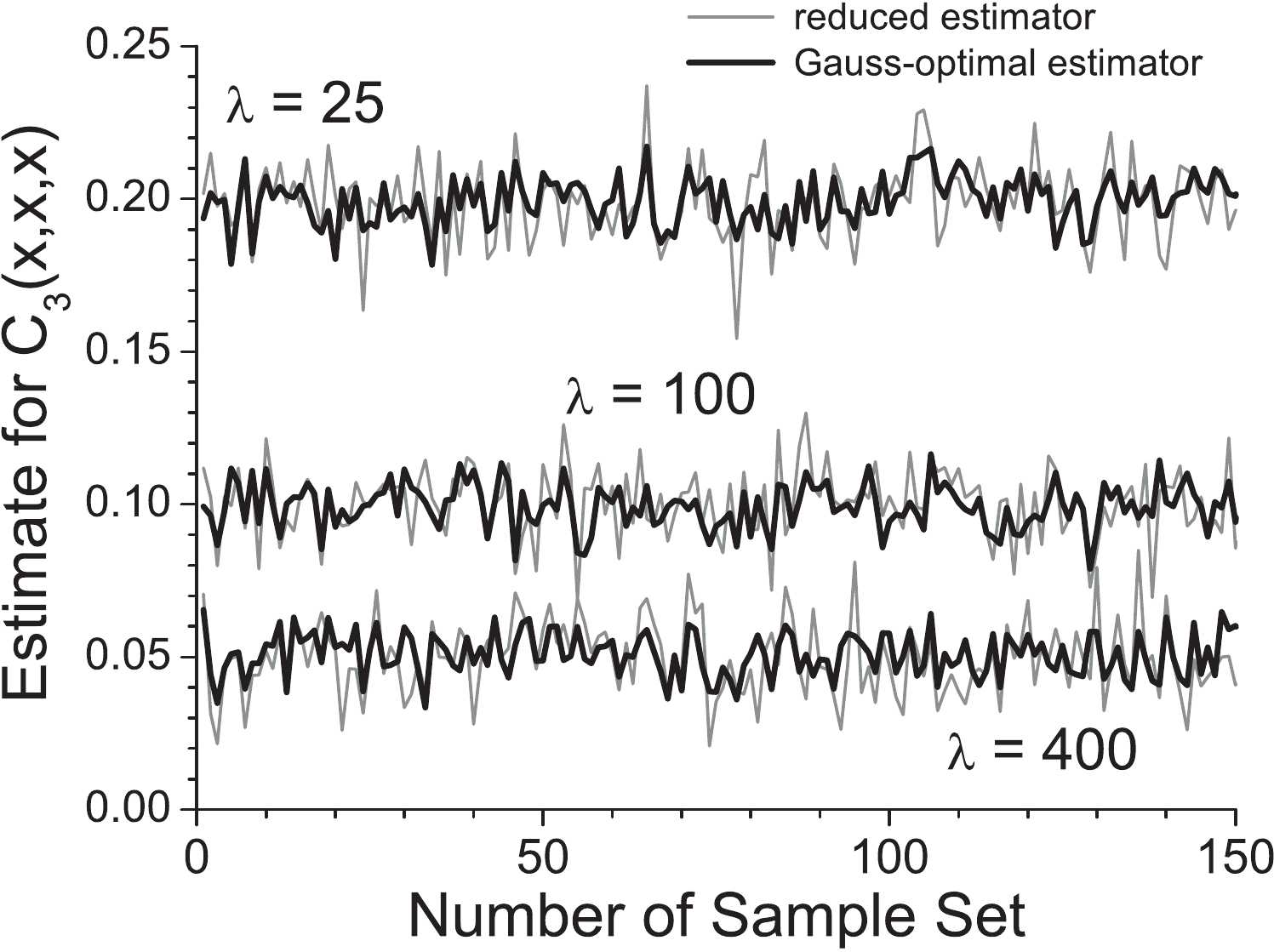}
	\caption{Estimates of $C_3(x,x,x) = \langle x^3 \rangle$ for three different quasi-Poisson distributed variables $x$
	 for increasing parameter $\lambda = 25,\,100,\,400$ (see text).  The new Gauss-optimal estimator  $c_3^{\rm (c,Go)}$ (black line) gives rise to less noise than the reduced estimator $c_3^{\rm (c)}$ (grey line).}
	\label{fig:Poisson}
\end{figure}
In the section above we showed that the estimator $c_3^{\rm (c)} = \overline{x^3}$ is not Gauss-optimal for estimating $C_3(x,x,x)$ of 
an average free random variable $x$. Instead we introduced the new estimator $c_3^{\rm (c,Go)}$ [see (\ref{eq:c3c_Go})].
 Here we compare the performance of 
 $c_3^{\rm (c)}$ and $c_3^{\rm (c,Go)}$ in a numerical experiment for a random variable $x$ that is derived from a 
 Poisson distributed variable $h$ with parameter $\lambda$ which possesses the well-known probability distribution
  $p_h(n,\lambda) = \lambda^n e^{-\lambda}/n!$ for finding the value $n$. A Poisson distribution is centered around $h = \lambda$ with a variance $\lambda$ and exhibits a slight asymmetry (skewness) that gives rise to a non-vanishing $\langle h^3 \rangle = \lambda^{1/2}$. 
The random variable 
\begin{equation}
  x = \lambda^{-1/2}(h - \lambda)
\end{equation}
   is constructed from $h$ in a way to insure $\langle x \rangle = 0 $. Moreover, the relations
\begin{eqnarray}
C_2(x,x) & = & \langle x^2 \rangle  =  1 \nonumber \\
C_3(x,x,x) & = & \langle x^3 \rangle  =  \lambda^{-1/2}.
\end{eqnarray}
hold as can be easily derived from the properties of the Poisson distributed variable $h$. 
Figure \ref{fig:QuasiPoisson} shows the probability distribution of $x$ for $\lambda = 25$ (where the discrete values for $x$ were replaced by vertical bars)
 in direct comparison with the corresponding Gauss-distribution $\exp(-x^2/2)/\sqrt{2\pi}$. The long right hand tail of the distribution clearly reveals an asymmetry that will lead to a non-vanishing $C_3(x,x,x)$. 
Figure \ref{fig:Poisson} shows 150 datapoints for the estimators $c_3^{\rm (c)}$ and $c_3^{\rm (c,Go)}$ that 
 were evaluated for different sets of  $m=10^5$ samples of $x$ for three different values of $\lambda$.
  The datapoints for $\lambda = 25,\,100,\,400$ scatter around the expected values for $C_3(x,x,x) = 0.2,\,0.1,\,0.05$. 
  Intriguingly, the values of the new estimator $c_3^{\rm (c,Go)}$ scatter significantly less than those of $c_3^{\rm (c)}$. 
  We estimated the variance of the datapoints from 300 samples and found good agreement with the 
  theoretically expected values of $\sigma^2 = m V(c_3^{\rm (c)})/10^5 \approx 1.5\times10^{-4}$ and
$\sigma^2 = m V(c_3^{\rm (c,Go)})/10^5 \approx 6\times10^{-5}$. 
  Overall the improvement of the variance of the scatter is a factor of $2.5$ in agreement with theory. 
  We emphasize the importance of this result for actual experiments.
   A factor of 2.5 less samples are sufficient to obtain the same certainty about an estimate of $C_3(x,x,x)$ using the new Gauss optimal estimator as compared to the reduced estimator.

\section{Consistency}
\label{sec:consistency}
Consistency of an estimator $c$ is its property to converge to the cumulant $C$ that it is estimating for $m \rightarrow \infty$. 
It is sufficient to show that the variance $V(c) \rightarrow 0$ for $m \rightarrow \infty$ \cite{blagouchineJSP2012}. We therefore find immediately consistency for all our estimators that are labeled with a V in Table \ref{OverviewEstimator}. Their variances tend to zero for increasing $m$ as can be seen from the expressions given in Appendix \ref{app:var}.
Consistency for any of our unbiased estimators $c$ is established by showing that their variances have the 
property
\begin{eqnarray}
 V(c) & =  & \langle c^2 \rangle - C^2 \nonumber \\
  &  = & \langle c^2 \rangle - \langle c \rangle^2 = O(m^{-1}), 
\end{eqnarray}
where $O(m^{-1})$ means that the order of all terms on the RHS is $O(m^{-1})$ or higher and consequently $V(c) \rightarrow 0$ for $m \rightarrow \infty$. 

We first note that all estimators are of the form
\begin{eqnarray}
 c & = & \sum_k c_k \nonumber \\
 c_k & = & \overline{p_{k,1}}\,\overline{p_{k,2}}\,...\overline{p_{k,n_k}}, \label{eq:ck}
\end{eqnarray} 
where the $p$s are polynomials of the random variables and $n_k$ is the number of factors that appear in $c_k$.
The statistical mean of $c_k$ is
\begin{equation}
 \langle c_k \rangle = \langle p_{k,1} \rangle \langle p_{k,2} \rangle ... \langle p_{k,n_k} \rangle + O(m^{-1}).  \label{c_leadingOrder}
\end{equation}
 The example $c_k = \overline{x} \ \overline{y} \ \overline{z} \ \overline{w}$ appears in (\ref{eq:c_example}) where
  the leading zero order contribution $\langle x \rangle \langle y \rangle \langle z \rangle \langle w \rangle$ is found.
The general relation (\ref{c_leadingOrder}) follows from considering (\ref{recursive}) where the leading order originates from the contributions for $\nu = 0$ in the (recursive) sum. We consequently find
\begin{equation}
 \langle c \rangle^2 = \sum_{k,k'} \langle p_{k,1} \rangle  ... \langle p_{k,n_k} \rangle \langle p_{k',1} \rangle
 ... \langle p_{k',n_{k'}} \rangle + O(m^{-1}). 
\end{equation}
 Since $c^2$ is also of the form defined in (\ref{eq:ck}) we find using again (\ref{c_leadingOrder})
 \begin{equation}
 \langle c^2 \rangle = \sum_{k,k'} \langle p_{k,1} \rangle  ... \langle p_{k,n_k} \rangle \langle p_{k',1} \rangle
 ... \langle p_{k',n_{k'}} \rangle + O(m^{-1}). 
\end{equation}
which with the result for $ \langle c \rangle^2$ establishes our claim $V(c) = \langle c^2 \rangle - \langle c \rangle^2 = O(m^{-1})$ and therefore consistency of any $c$.
\section{Conclusion}
In conclusion, we derived multivariate unbiased estimators for the second-, third- and fourth-order cumulants
 including several cases with average-free variables and pairs of variables with vanishing second order cumulant. 
 The reduced third order estimators  $c_3^{\rm (c)} = \overline{x^3}$ and $c_3^{\rm (h)} = \overline{x^2 z}$ for average free
  variables turned surprisingly out to be not Gauss-optimal, while Gauss-optimal alternative estimators could be derived in Section \ref{sec:go}.  An overview over the estimators along 
   with a new nomenclature is given in Table \ref{OverviewEstimator}. 
   As a side result of our work we gave two simple recursive formulas for finding multivariate cumulants from moments and vice versa.
   We expect that some of the new estimators will soon find application in signal processing especially for estimating higher order noise spectra from cumulants of Fourier coefficients.

\appendices
\section*{Acknowledgment}
The authors acknowledge financial support
of the Deutsche Forschungsgemeinschaft under Grant No. HA
3003/7-1.

\section{Variances of Estimators}
\label{app:var}
Here we express the variances $V(c)$ of a selection of estimators $c$ in terms of higher order cumulants. 
The calculation of $V(c) = \langle c^2 \rangle - \langle c \rangle^2$ requires the evaluation of 
many terms like $\langle \overline{xy}\, \overline{y}\,\overline{xz}\overline{z}\rangle$ etc. similar
 to terms in Section \ref{sec:c3} and \ref{sec:c4} where 
 such expression had been 
 evaluated by hand in terms of moments. We employ here Computer Algebra to first express 
 the expected statistical values in terms of moments. The method in Section \ref{sec:c4} suggests a recursive algorithm.
We seek to calculate
\begin{equation}
   \langle \overline{p_1} \, \overline{p_2} \cdots  \overline{p_N} \rangle   
\end{equation}
where $p_j$ are polynomials of the stochastic variables $x$, $y$, etc.
The recursive method requires us to be specific about the number $m$ of samples that is used for calculating $\overline{p_j}$. We therefore introduce the notation $A_m(p_j) = \overline{p_j}$ which keeps track of $m$.
The calculation of 
   $\langle A_m(p_1) \cdots A_m(p_N) \rangle$ can be reformulated with the recursive helper function
\begin{eqnarray}
 H(A_m(p_1), \cdots, A_m(p_N)) &  & \nonumber \\
 &   & \hspace{-3cm}=\sum_{\nu = 0}^{N-1} \hspace{-2cm}  \sum_{\tiny \hspace{2cm}\begin{array}{ll}\textrm{all partitions of}\,2,3,...,N\, \textrm{into} \nonumber \\
 ((i_1,...,i_\nu),(j_1,...,j_{N-1-\nu}))\end{array}} \hspace{-2cm} m \langle p_1 p_{i_1} \cdots p_{i_\nu}\rangle \nonumber  \\
 & & \hspace{-2cm}\times H(A_{m-1}(p_{j_1}), \cdots, A_{m-1}(p_{j_{n-1 - \nu}})) \label{recursive}
\end{eqnarray}
and
\begin{eqnarray}
\langle \overline{p_1} \, \overline{p_2} \cdots  \overline{p_N} \rangle   = \frac{1}{m^N}H(A_m(p_1), \cdots, A_m(p_N)). \label{evalEst}
\end{eqnarray}    
The partitions under the sum means that the indices $2$ to $N$ have to be partitioned into $\nu$ different indices $i$ and into $N-\nu-1$ indices $j$.  There are ${\tiny \left( \begin{array}{c}   N-1 \\ \nu \end{array}\right)}$ different possible partitions. 
All possible partitions $\{i_1,...,i_\nu\}$ of a set $A =\{1,2,...,N\}$ can be obtained in the computer algebra system MATHEMATICA via the function ${\rm Subsets}[A,\nu]$.
 The factor $m$ in (\ref{recursive}) leads after recursion to factors of the form $m(m-1) \cdots$ [compare e.g. (\ref{EqExpAv})]. For $\nu = 0$ the first factor in the sum is $\langle p_1 \rangle$. For $\nu = N-1$ there are no indices $j$. For that case $H() = 1$ has to be defined.
After $V(c)$ is expressed in terms of higher order moments with the help of (\ref{evalEst}) the moments are expressed in terms of higher order cumulants using another recursive method 
\begin{eqnarray}
M_N(x_1,\cdots,x_N) = \hspace{1cm}    &   & \hspace{-1.1cm}  C_N(x_1,...,x_N)  \nonumber \\
 &   & \hspace{-4.5cm} +\sum_{\nu = 1}^{N-1} \hspace{-2cm} \sum_{\tiny \hspace{2cm}\begin{array}{ll}\textrm{all partitions of}\,1,2,...,N\, \textrm{into} \\ 
 ((i_1,...,i_\nu),(j_1,...,j_{N-\nu}))\end{array}} \hspace{-2.5cm} \frac{\nu}{N}  C_{\nu}(x_{i_1},...,x_{i_\nu})  M_{N-\nu}(x_{j_1},\cdots x_{j_{N-\nu}}) \label{Cum}
\end{eqnarray}
(a short proof is given in Appendix \ref{app:rec}).

For completeness, we also state the inverse formula (cumulant generating formula)
\begin{eqnarray}
C_N(x_1,\cdots,x_N) = \hspace{1cm}    &   & \hspace{-1.1cm}  M_N(x_1,\cdots,x_N)  \nonumber \\
 &   & \hspace{-4.5cm} -\sum_{\nu = 1}^{N-1} \hspace{-2cm} \sum_{\tiny \hspace{2cm}\begin{array}{ll}\textrm{all partitions of}\,1,2,...,N\, \textrm{into} \\
 ((i_1,...,i_\nu),(j_1,...,j_{N-\nu}))\end{array}} \hspace{-2.5cm} \frac{\nu}{N}  C_{\nu}(x_{i_1},...,x_{i_\nu})  M_{N-\nu}(x_{j_1},...,x_{j_{N-\nu}}).\label{Mom}
\end{eqnarray}
Another recursive method for obtaining  multivariate $M_N(x_1,\cdots,x_N)$ and $C_N(x_1,\cdots,x_N)$ 
had been given by Smith before \cite{smithTAS1995}. He required a formula with multiple sums over indices  for the case of moments
(instead of only two sums in our case). The case of cumulants required in addition a combination of two similarly complex formulas. We are not aware that our more simple form had been given in literature before. A non-recursive method due to Leonov and Shiryaev for calculating multivariate cumulants or moments using multiple partitions (which albeit may need to be constructed recursively) can be found in \cite{leonovTPA1959} and \cite{mendelIEEE1991}.
 
After a computer algebra implementation of (\ref{evalEst}) and (\ref{Cum}) we found the variances below. 
The results for the k-statistics $V(c_3^{\rm (b)})$ and $V(c_4^{\rm (b)})$ are in agreement with \cite{weissteinMATHEMATICA}.
\begin{equation}
 V(c_2^{\rm (a)}) = \frac{C_4(x,x,y,y)}{m}  + \frac{C^2_2(x,y)}{m-1}+\frac{C_2(x,x)C_2(y,y)}{m-1}
\end{equation}
\begin{eqnarray}
 V(c_2^{\rm (d)}) & = & \langle x^2 y^2 \rangle /m - \langle x y\rangle^2 / m \nonumber \\
    &  = & \frac{C_4(x,x,y,y)}{m}  + \frac{C^2_2(x,y)}{m}+\frac{C_2(x,x)C_2(y,y)}{m} \nonumber \\
    & &
\end{eqnarray}
\begin{eqnarray}
 V(c_3^{\rm (a)}) & = & \frac{C_6(x,x,y,y,z,z)}{m} \nonumber \\
 & & + \left( \frac{C_4(x,x,y,y)C_2(z,z)}{m-1} + \textrm{2 o.p.} \right) \nonumber \\
 & & + \left( \frac{2 C_4(x,x,y,z)C_2(y,z)}{m-1} + \textrm{2 o.p.} \right) \nonumber \\
 & & + \left( \frac{2 C_3(x,x,y)C_3(y,z,z)}{m-1} + \textrm{2 o.p.} \right) \nonumber \\
 & & + \frac{3 C^2_3(x,y,z)}{m-1} \nonumber \\
 & & + \left(\frac{m C_2(x,x)C^2_2(y,z)}{(m-1)(m-2)} + \textrm{2 o.p.}\right)\ \nonumber \\
  & & + \frac{2 m C_2(x,y)C_2(y,z)C_2(z,x)}{(m-1)(m-2)} \nonumber \\
 & & + \frac{m C_2(x,x)C_2(y,y)C_2(z,z)}{(m-1)(m-2)}
  \end{eqnarray}
\begin{eqnarray}
 V(c_3^{\rm (b)}) & = & \frac{C_6(x,x,x,x,x,x)}{m} \nonumber \\
 & & +  \frac{9 C_4(x,x,x,x)C_2(x,x)}{m-1}  \nonumber \\
 & & + \frac{9 C^2_3(x,x,x)}{m-1} \nonumber \\
 & & + \frac{6 m C^3_2(x,x)}{(m-1)(m-2)}
  \end{eqnarray}  
\begin{eqnarray}
 V(c_3^{\rm (c)}) & = & (\langle x^6\rangle -  \langle x^3\rangle^2)/m \nonumber \\
 & = & \frac{C_6(x,x,x,x,x,x)}{m} \nonumber \\
 & & +  \frac{15 C_4(x,x,x,x)C_2(x,x)}{m} \nonumber \\
 & & + \frac{9 C^2_3(x,x,x)}{m}  + \frac{15 C^3_2(x,x)}{m}
\end{eqnarray}   
\begin{eqnarray}
 V(c_3^{\rm (d)}) & = & \frac{C_6(x,x,y,y,z,z)}{m} \nonumber \\
 & & + \left( \frac{C_4(x,x,y,y)C_2(z,z)}{m} + \textrm{2 o.p.} \right) \nonumber \\
 & & + \left( \frac{4 C_4(x,x,y,z)C_2(y,z)}{m} + \textrm{2 o.p.} \right) \nonumber \\
 & & + \left( \frac{2 C_3(x,x,y)C_3(y,z,z)}{m} + \textrm{2 o.p.} \right) \nonumber \\
 & & + \frac{3 C^2_3(x,y,z)}{m} \nonumber \\
 & & + \left(2\frac{ C_2(x,x)C^2_2(y,z)}{m} + \textrm{2 o.p.}\right)\ \nonumber \\
  & & + \frac{8  C_2(x,y)C_2(y,z)C_2(z,x)}{m} \nonumber \\
 & & + \frac{C_2(x,x)C_2(y,y)C_2(z,z)}{m}
  \end{eqnarray}  
\begin{eqnarray}
 V(c_4^{(b)})  & = & \frac{C_8(x,x,x,x,x,x,x,x)}{m} \nonumber \\
 & & +  \frac{16 C_6(x,x,x,x,x,x)C_2(x,x)}{m-1} \nonumber \\
  & & +  \frac{48 C_5(x,x,x,x,x)C_3(x,x,x)}{m-1} \nonumber \\
    & & + \frac{34  C^2_4(x,x,x,x) }{m-1} \nonumber \\
 & & + \frac{72 m C_4(x,x,x,x) C^2_2(x,x)}{(m-1)(m-2)} \nonumber \\
 & & + \frac{144 m C^2_3(x,x,x) C_2(x,x)}{(m-1)(m-2)} \nonumber \\ 
 & & + \frac{24 m (m+1) C^4_2(x,x)}{(m-1)(m-2)(m-3)}  
  \end{eqnarray}  
The variance of the  unbiased estimator $c_4^{(c)}(x)$ where $\langle x \rangle = 0$ 
was given by Blagouchine in \cite{BlagouchineIEEE2009} in terms of moments of $x$. We could verify their result and 
rewrite it here in terms of cumulants
\begin{eqnarray}
 V(c_4^{(c)})  & = & \frac{C_8(x,x,x,x,x,x,x,x)}{m} \nonumber \\
 & & +  \frac{16 C_6(x,x,x,x,x,x)C_2(x,x)}{m} \nonumber \\
  & & +  \frac{56 C_5(x,x,x,x,x)C_3(x,x,x)}{m} \nonumber \\
    & & + \frac{(34 m -16) C^2_4(x,x,x,x) }{m(m-1)} \nonumber \\
 & & + \frac{72 C_4(x,x,x,x) C^2_2(x,x)}{m-1} \nonumber \\
 & & + \frac{160 C^2_3(x,x,x) C_2(x,x)}{m} \nonumber \\ 
 & & + \frac{24 (m+2) C^4_2(x,x)}{m(m-1)}  
  \end{eqnarray}

  \begin{eqnarray}
 V(c_4^{(ca)})  & = & \frac{C_8(a,a^*,a,a^*,b,b^*,b,b^*)}{m} \nonumber \\
 & & +  \frac{2 C_6(a,a^*,b,b^*,b,b^*)C_2(a,a^*)}{m} \nonumber \\
  & & +  \frac{2 C_6(a,a^*,a,a^*,b,b^*)C_2(b,b^*)}{m} \nonumber \\
      & & + \frac{C_4(a,a^*,a,a^*)C_4(b,b^*,b,b^*) }{m-1} \nonumber \\
            & & + \frac{(7m-6) C^2_4(a,a^*,b,b^*) }{m(m-1)} \nonumber \\
    & & + \frac{C_4(a,a^*,a,a^*)C^2_2(b,b^*) }{m-1} \nonumber \\
        & & + \frac{4 C_4(a,a^*,b,b^*)C_2(a,a^*)C_2(b,b^*) }{m} \nonumber \\
                & & + \frac{ C_4(b,b^*,b,b^*)C^2_2(a,a^*) }{m-1} \nonumber \\
    & & + \frac{ C^2_2(a,a^*)C^2_2(b,b^*)}{m-1} 
  \end{eqnarray}  
    \begin{eqnarray}
 V(c_4^{(cb)})  & = & \cdots \nonumber \\
 & + &  \frac{ (m+1) C^2_2(a,a^*)C^2_2(b,b^*)}{m(m-1)}\nonumber \\
  & + &  \frac{ 2(m+1) C_2(a,a^*)C_2(a,b^*) C_2(b,a^*)C_2(b,b^*)}{m(m-1)}\nonumber \\
    & & + \frac{ (m+1)C^2_2(a,b^*)C^2_2(b,a^*)}{m(m-1)} 
  \end{eqnarray}
In the last result we omitted the quite lengthy contributions of terms that included cumulants of order three and higher.   
\section{Recursive Calculation of Multivariate Cumulants}
\label{app:rec}
Here we give a short proof of the cumulant/moment generating recursive formulas used in Appendix \ref{app:var}.
Smith gives the following recursive formula for the univariate case \cite{smithTAS1995}
\begin{eqnarray}
C_N(u,...,u) & = & M_N(u,...,u) \nonumber \\
 & & \hspace{-2cm} - \sum_{\nu = 1}^{N-1}{\small \left( \begin{array}{c} N-1 \\ \nu \end{array} \right)}C_{N-\nu}(u,...,u)M_{\nu}(u,...,u)
\end{eqnarray}
or equivalently 
\begin{eqnarray}
C_N(u,...,u) & = & M_N(u,...,u) \nonumber \\
 & & \hspace{-2cm} - \sum_{\nu = 1}^{N-1}{\small \left( \begin{array}{c} N-1 \\ \nu-1 \end{array} \right)}C_{\nu}(u,...,u)M_{N-\nu}(u,...,u).
 \label{rekWiki}
\end{eqnarray}
The multivariate case is obtained from (\ref{rekWiki}) considering $u = \vec{k}\vec{x}$.
The coefficient of the term with the factor $k_1 k_2 ... k_N$ in $C_N( \vec{k}\vec{x},..., \vec{k}\vec{x})$ [LHS of (\ref{rekWiki})] is
$N!C_N(x_1,...,x_N)$ where we made use of the multilinearity of cumulants. 
The coefficient of the RHS of (\ref{rekWiki}) is
\begin{eqnarray}
   &   & N!  M_N(x_1,...,x_N)  \nonumber \\
 &   &  -\sum_{\nu = 1}^{N-1} \hspace{-2cm} \sum_{\tiny \hspace{2cm}\begin{array}{ll}\textrm{all permutations of}\,1,2,...,N\, \textrm{into} \nonumber \\
 ((i_1,...,i_\nu),(j_1,...,j_{N-\nu}))\end{array}} \hspace{-2cm} {\small \left( \begin{array}{c} N-1 \\ \nu-1 \end{array} \right)}   \nonumber \\
 & &\hspace{1cm}\times C_\nu(x_{i_1},..., x_{i_\nu})  M_{N-\nu}(x_{j_1},...,x_{j_{N-\nu}}).
\end{eqnarray}
Since $C_\nu(x_{i_1},..., x_{i_\nu})$ and $M_{N-\nu}(x_{j_1},...,x_{j_{N-\nu}})$ are identical under permutation of their arguments,
we can rewrite the above equation as 
\begin{eqnarray}
   &   & N!  M_N(x_1,...,x_N)  \nonumber \\
 &   &  -\sum_{\nu = 1}^{N-1} \hspace{-2cm} \sum_{\tiny \hspace{2cm}\begin{array}{ll}\textrm{all partitions of}\,1,2,...,N\, \textrm{into} \nonumber \\
 ((i_1,...,i_\nu),(j_1,...,j_{N-\nu}))\end{array}} \hspace{-2cm} {\small \left( \begin{array}{c} N-1 \\ \nu-1 \end{array} \right)} \nu! (N-\nu)!  \nonumber \\
 & &\hspace{1cm}\times C_\nu(x_{i_1},..., x_{i_\nu})  M_{N-\nu}(x_{j_1},...,x_{j_{N-\nu}}),
\end{eqnarray}
where the factors $\nu!$ and $(N-\nu)!$ correctly regard the multiplicities of identical factors $C_\nu$ and $M_{N-\nu}$. 
After rewriting ${\small \left( \begin{array}{c} N-1 \\ \nu-1 \end{array} \right)} = (N-1)!/((\nu-1)!(N-\nu)!)$ we obtain (\ref{Mom}).
The moment generating formula (\ref{Cum}) follows directly from  (\ref{Mom}) by rearranging the sums from the RHS to the LHS.

\ifCLASSOPTIONcaptionsoff
  \newpage
\fi


\bibliographystyle{IEEEtran}
\bibliography{references}

\begin{thebibliography}{10}
\providecommand{\url}[1]{#1}
\csname url@samestyle\endcsname
\providecommand{\newblock}{\relax}
\providecommand{\bibinfo}[2]{#2}
\providecommand{\BIBentrySTDinterwordspacing}{\spaceskip=0pt\relax}
\providecommand{\BIBentryALTinterwordstretchfactor}{4}
\providecommand{\BIBentryALTinterwordspacing}{\spaceskip=\fontdimen2\font plus
\BIBentryALTinterwordstretchfactor\fontdimen3\font minus
  \fontdimen4\font\relax}
\providecommand{\BIBforeignlanguage}[2]{{%
\expandafter\ifx\csname l@#1\endcsname\relax
\typeout{** WARNING: IEEEtran.bst: No hyphenation pattern has been}%
\typeout{** loaded for the language `#1'. Using the pattern for}%
\typeout{** the default language instead.}%
\else
\language=\csname l@#1\endcsname
\fi
#2}}
\providecommand{\BIBdecl}{\relax}
\BIBdecl

\bibitem{comonSP1994}
P.~Comon, ``{Independent component analysis, A new concept?}'' \emph{Signal
  Process.}, vol.~36, p. 287, 1994.

\bibitem{brillingerAMS1965}
D.~R. Brillinger, ``An introduction to polyspectra,'' \emph{Ann. Math.
  Statist.}, vol.~36, p. 1351, 1965.

\bibitem{hagelePRB2018}
D.~H\"agele and F.~Schefczik, ``Higher-order moments, cumulants, and spectra of
  continuous quantum noise measurements,'' \emph{Phys. Rev. B}, vol.~98, p.
  205143, 2018.

\bibitem{FisherPLMS1928}
R.~A. Fisher, ``Moments and product moments of sampling distributions,''
  \emph{Proceedings of the London Mathematical Society}, vol. s2-30, no.~1, pp.
  199--238, 1928.

\bibitem{DiNardoBernoulli2008}
\BIBentryALTinterwordspacing
E.~Di~Nardo, G.~Guarino, and D.~Senato, ``A unifying framework for k
  -statistics, polykays and their multivariate generalizations,''
  \emph{Bernoulli}, vol.~14, no.~2, pp. 440--468, 05 2008. [Online]. Available:
  \url{https://doi.org/10.3150/07-BEJ6163}
\BIBentrySTDinterwordspacing

\bibitem{rotaSIAM1994}
G.-C. Rota and B.~D. Taylor, ``{The Classical Umbral Calculus},'' \emph{SIAM J.
  Math. Anal.}, vol.~25, p. 694, 1994.

\bibitem{mansourCONF1998}
A.~Mansour, A.~K. Kardec~Barros, and N.~Ohnishi, ``Comparison among three
  estimators for high order statistics,'' \emph{Fifth International Conference
  on Neural Information}, p. 899, 1998.

\bibitem{BlagouchineIEEE2009}
I.~V. Blagouchine and E.~Moreau, ``Unbiased adaptive estimations of the
  fourth-order cumulant for real random zero-mean signal,'' \emph{IEEE
  Transactions on Signal Processing}, vol.~57, no.~9, pp. 3330--3346, Sept
  2009.

\bibitem{gardinerBOOK2009}
C.~Gardiner, \emph{Stochastic Methods}, 4th~ed.\hskip 1em plus 0.5em minus
  0.4em\relax Berlin Heidelberg: Springer, 2009.

\bibitem{starosielecBOOK2012}
S.~Starosielec, \emph{Rauschspektroskopie h\"oherer Ordnungen
  (Dissertation)}.\hskip 1em plus 0.5em minus 0.4em\relax Sierke Verlag, 2012.

\bibitem{starosielecRSI2010}
S.~Starosielec, R.~Fainblat, J.~Rudolph, and D.~H\"agele, ``Two-dimensional
  higher order noise spectroscopy up to radio frequencies,'' \emph{Rev.
  Scientific Instrum.}, vol.~81, p. 125101, 2010.

\bibitem{gabelliNJP2013}
J.~Gabelli, L.~Spietz, J.~Aumentado, and B.~Reulet, ``Electron–photon
  correlations and the third moment of quantum noise,'' \emph{New J. Phys.},
  vol.~15, p. 113045, 2013.

\bibitem{blagouchineJSP2012}
I.~V. Blagouchine and E.~Moreau, ``{Comments on 'Unbiased estimates for moments
  and cumulants in linear regression'},'' \emph{J. Stat. Plan. Infer.}, vol.
  142, p. 1027, 2012.

\bibitem{smithTAS1995}
P.~J. Smith, ``A recursive formulation of the old problem of obtaining moments
  from cumulants and vice versa,'' \emph{Amer. Statistician}, vol.~49, p. 217,
  1995.

\bibitem{leonovTPA1959}
V.~P. Leonov and A.~N. Shiryaev, ``{On a Method of Calculation of
  Semi-Invariants},'' \emph{Theory Probab. Appl.}, vol.~4, p. 319, 1959.

\bibitem{mendelIEEE1991}
J.~M. Mendel, ``{Tutorial on Higher-Order Statistics (Spectra) in Signal
  Processing and System Theory: Theoretical Results and Some Applications},''
  \emph{Proc. IEEE}, vol.~19, p. 278, 1991.

\bibitem{weissteinMATHEMATICA}
{Weisstein, Eric W. 'k-Statistic.' From MathWorld--A Wolfram Web Resource.
  http://mathworld.wolfram.com/k-Statistic.html}.

\end{thebibliography}

\begin{IEEEbiography}[{\includegraphics[width=1in,height=1.25in,clip,keepaspectratio]{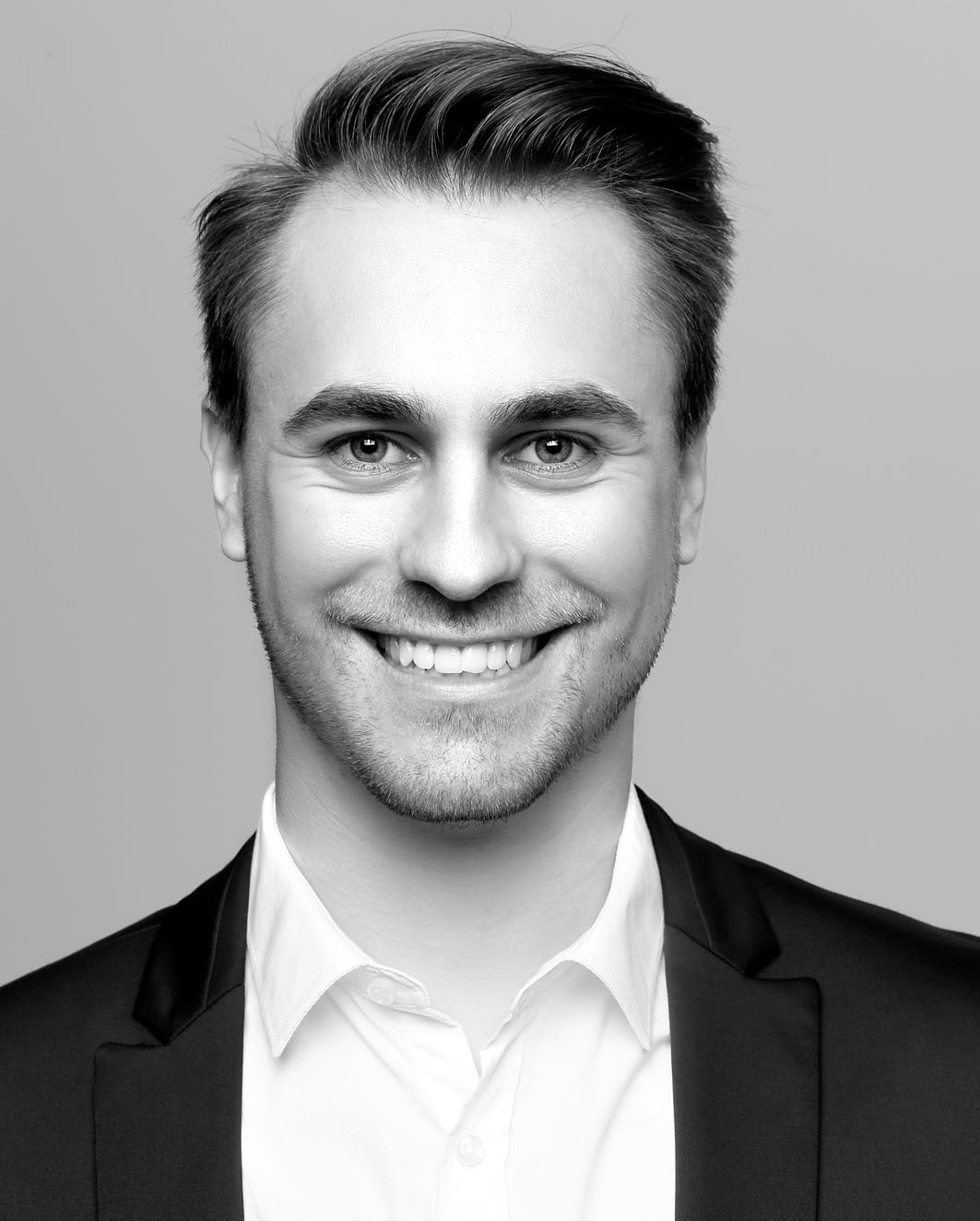}}]{Fabian Schefczik}
studied physics (B.Sc. in 2014 and M.Sc. in 2016) with a focus on hadron physics at Ruhr-University Bochum (Germany). He currently pursues a doctorate in the research group of D. H\"agele on a topic including quantum measurements and higher order noise spectroscopy.\end{IEEEbiography}

\begin{IEEEbiography}[{\includegraphics[width=1in,height=1.25in,clip,keepaspectratio]{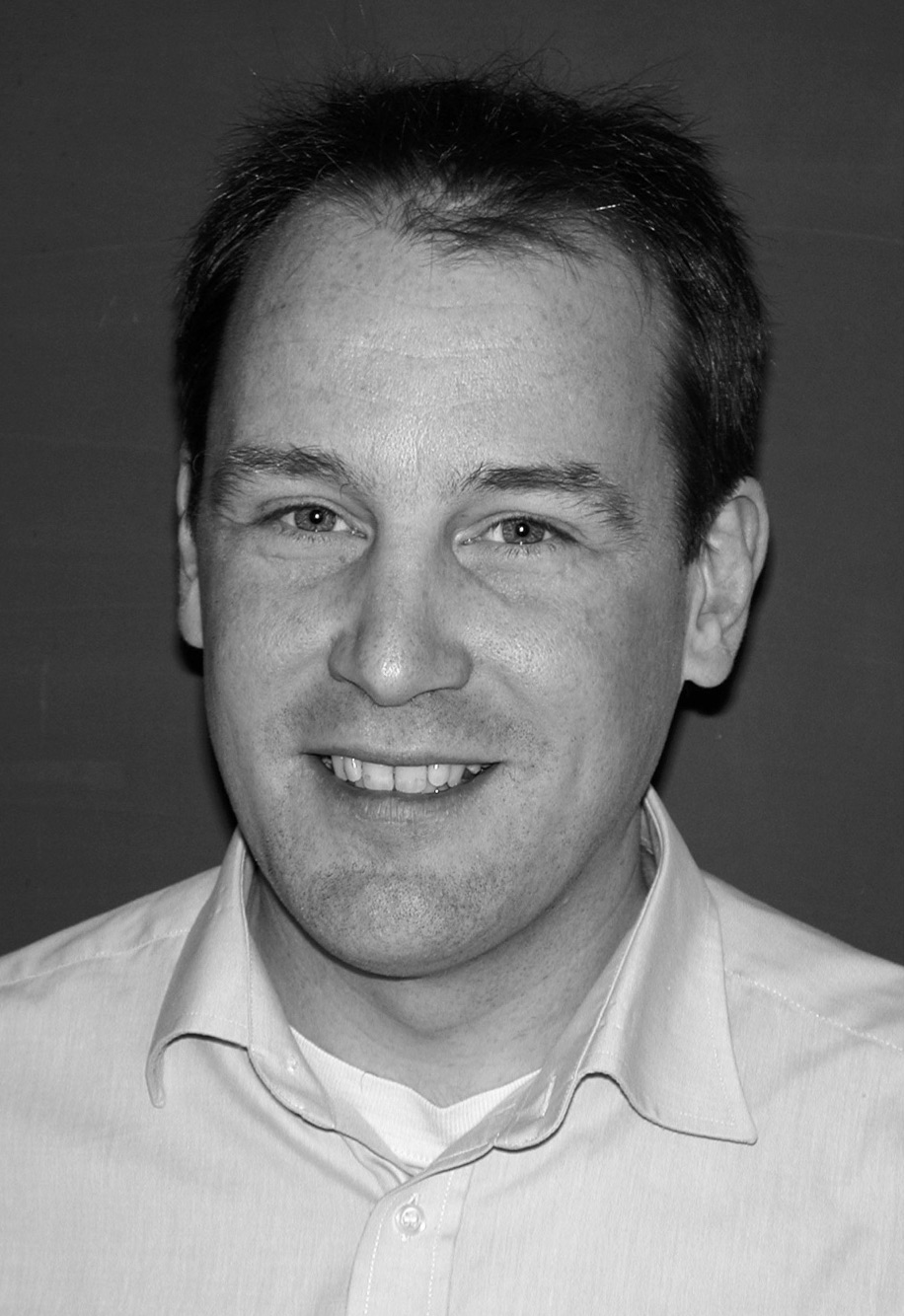}}]{Daniel H\"agele} received his Diploma in physics from the University of Ulm (Germany) in 1995 and his  PhD in physics form the Philipps University Marburg (Germany) in 1999. In between he spent 13 months in alternative civilian service. After joining Lawrence Berkeley Laboratory (California) as a guest scientist in 2000 he became a PostDoc at Leibnitz University Hannover (Germany) in 2002. Since 2006 he is a professor for spectroscopy of condensed matter in the Department of Physics and Astronomy at the Ruhr University Bochum (Germany). 
His research interests include solid state calorics, semiconductor spin physics, continuous quantum measurements, and real-time measurements of higher order spectra with GHz bandwidth. 
\end{IEEEbiography}

\end{document}